\documentclass[12pt]{article}

\usepackage{epsfig}
\usepackage{amsmath}
\usepackage{amsfonts}
\usepackage{mathbbol}
\usepackage{amssymb}
\usepackage{xcolor}
\usepackage{mathrsfs}
\usepackage{hyperref}

\newtheorem{thm}{Theorem}[section]

\newtheorem{lem}[thm]{Lemma}

\newtheorem{ass}[thm]{Assumption}
\newtheorem{de}[thm]{Definition}

\newtheorem{rem}[thm]{Remark}
\newtheorem{problem}[thm]{Problem}

\numberwithin{equation}{section}

\newcommand{\be}{\begin{equation*}}
\newcommand{\ben}{\begin{equation}}
\newcommand{\ee}{\end{equation*}}
\newcommand{\een}{\end{equation}}

\newcommand{\esssup}{\operatorname{ess\,sup}}
\newcommand{\cc}{{\mathrm c}}

\DeclareSymbolFontAlphabet{\amsmathbb}{AMSb}

\newcommand{\bbr}{\mathbb{R}}

\newcommand{\hsym}{{\varphi}}
\newcommand{\usym}{{\psi}}
\newcommand{\gbG}{\Phi}
\newcommand{\bgbG}{\mathbb{\Phi}}
\newcommand{\Usym}{\Psi}
\newcommand{\bUsym}{\mathbb{\Psi}}
\newcommand{\hfn}{{h}}

\newcommand{\EXP}{\operatorname{\textsf{\upshape E}}}
\newcommand{\PR}{\operatorname{\textsf{\upshape P}}}

\newcommand{\inv}{{\rm inv}}

\newcommand{\acal}{\mathcal{A}}
\newcommand{\bcal}{\mathcal{B}}
\newcommand{\ccal}{\mathcal{C}}
\newcommand{\fcal}{\mathcal{F}}

\newcommand{\ical}{\mathcal{I}}

\newcommand{\pcal}{\mathcal{P}}
\newcommand{\scal}{\mathcal{S}}

\newcommand{\vcal}{\mathcal{V}}

\newcommand{\xcal}{\mathcal{X}}

\newcommand{\iscr}{\mathscr I}
\newcommand{\lscr}{\mathscr L}

\newcommand{\di}{\mathrm{d}}

\begin{document}

\title{A Dynamic Competitive Equilibrium Model of
Irreversible Capacity Investment with Stochastic
Demand and Heterogeneous Producers}
\author{Constantinos Kardaras,
Alexandros Pavlis and Mihail Zervos}
\date{\today}
\maketitle

\begin{abstract} \noindent
We formulate a continuous-time competitive equilibrium
model of irreversible capacity investment in which a
continuum of heterogeneous producers supplies a single
non-durable good subject to exogenous stochastic
demand.
Each producer optimally adjusts both output and
capacity over time in response to endogenous price
signals, while investment decisions are irreversible.
Market clearing holds continuously, with prices
evolving endogenously to balance aggregate supply
and demand through a constant-elasticity demand
function driven by a stochastic base component.
The model admits a mean-field interpretation, as
each producer’s decisions both influence and are
influenced by the aggregate behaviour of all others.
We show that the equilibrium price process can be
expressed as a nonlinear functional of the exogenous
base demand, leading to a three-dimensional
singular stochastic control problem for each
producer.
We derive an explicit solution to the associated
Hamilton-Jacobi-Bellman equation, including a
closed-form characterisation of the free-boundary
surface separating investment and waiting regions.
\smallskip

\smallskip

\end{abstract}

\section{Introduction}

We consider a continuum of heterogeneous producers
of a single non-durable good, subject to exogenous
stochastic demand in continuous time.
Each producer determines, over time, both the level of
output and the amount of investment aimed at improving
production efficiency.
Producers must commit to their choices, meaning that
investments are irreversible.
This assumption is particularly plausible in industries
with high upfront costs, where the resale of capital may
entail substantial losses.

The price of the good is determined in equilibrium as the
aggregate outcome of individual producers' actions.
In particular, each producer acts as a price taker, meaning
that they have no direct influence on the price process.
In this competitive setting, the decision problem faced by
each producer gives rise to an infinite-horizon singular
stochastic control problem.

Continuous time irreversible capacity expansion models
featuring a single producer of a good with exogenous
prices or demand have attracted considerable attention
in the literature.
Early references in economics include the 
pioneering work of Manne~\cite{M61} (see Van
Mieghem~\cite{VM03} for a survey).
In the mathematical literature, models of stochastic
control with irreversible investment have been studied
by
Davis, Dempster, Sethi and Vermes~\cite{DDSV87},
Arntzen~\cite{A95},
{\O}ksendal~\cite{O00},
Wang~\cite{W03},
Chiarolla and Haussmann~\cite{CH05},
Alvarez~\cite{A10},
L{\o}kka and Zervos~\cite{LZ11},
Chiarolla and Ferrari~\cite{CF14},
Ferrari~\cite{F15},
De~Angelis, Federico and Ferrari~\cite{DAFF17},
Federico, Rosestolato and Tacconi~\cite{FRT19},
Koch and Vargiolu~\cite{KV21},
Dammann and Ferrari~\cite{DF22},
Han and Yi~\cite{HY22},
among others, listed here in rough alphabetical order
(see also references therein).
Furthermore, capacity expansion models allowing
costly reversibility were introduced by
Abel and Eberly~\cite{AE96}
and were subsequently studied by
Guo and Pham~\cite{GP05},
Merhi and Zervos~\cite{MZ07},
L{\o}kka and Zervos~\cite{LZ13},
De~Angelis and Ferrari~\cite{DAF14},
Federico and Pham~\cite{FP14}, and
Han, Yi and Zhang~\cite{HYZ24}.

Oligopolistic games of irreversible investment in
continuous time have been studied by
Back and Paulsen~\cite{BP09}
and by Steg~\cite{S12, S24}.
The models in these works consider settings in
which identical producers are exposed to an
exogenous economic shock modelled as a
one-dimensional It\^{o} diffusion.
Each individual producer solves a singular stochastic
control problem in which the running payoff depends
on both the exogenous shock and the capital stocks
of all market participants.

Recently, several mean-field investment games
have been analysed.
These include the finite-horizon irreversible
capacity expansion model studied by
Campi, De~Angelis, Ghio and Livieri~\cite{CDAGL22},
in which the drift of the price of the good produced by each
firm depends on the average installed capacity across
all producers.
Cao and Guo~\cite{CG22} study an infinite-horizon
partially reversible capacity expansion model
in which firms' production levels are modelled
as independent copies of a controlled geometric
Brownian motion, while the running payoff of each
producer depends on the distribution of the
long-run production level across all producers.
In a similar vein, A\"id, Basei and Ferrari~\cite{ABF25}
investigate an infinite-horizon irreversible capacity
expansion model in which firms' production levels
are modelled as independent copies of a controlled
geometric Brownian motion with regime switching.

To the best of our knowledge, the model we study is
the first in the capacity investment literature to
exhibit the key features of a classical dynamic
market equilibrium problem.
The demand for the good follows a constant-elasticity
demand function, where the base demand evolves
as a geometric Brownian motion and the prevailing price
determines the realised demand level.
In response to price signals, producers dynamically
adjust both their production capacities and output
rates over time.
Moreover, market clearing holds at all times:
prices evolve endogenously to maintain a continuous
balance between aggregate supply and demand.

Alternatively, the model can be interpreted as a
mean-field game, in which each producer's decisions
both influence and are influenced by the aggregate
behaviour of all other producers.
Unlike the capacity investment game models
discussed above, producers are not statistically
identical.
In particular, heterogeneity arises from
differences in initial production capacities, individual
discount rates, and unit costs of capacity expansion
as well as production costs.

We determine a Markovian market equilibrium by
first solving a nonlinear functional equation that
expresses the market-clearing price process of
the good as a function of the exogenous base
demand process.
This solution allows us to represent each producer's
objective in terms of a suitable geometric Brownian
motion and its running maximum.
Building on this result, we reformulate each producer's
optimisation problem as a classical three-dimensional
singular stochastic control problem.

We solve the singular stochastic control problem that
characterises each producer's optimal decisions in
equilibrium by first deriving a suitable explicit solution
to the associated Hamilton-Jacobi-Bellman (HJB)
equation.
In particular, we obtain a closed-form expression for
the free-boundary surface that separates the so-called
``investment'' and ``wait'' regions (see the surface
$\scal_i^1$, depicted by the purple lines in
Figure~\ref{pic1}).
To the best of our knowledge, this is the first
three-dimensional stochastic control problem with
a non-trivial explicit solution to have been studied
in the literature.

\section{Individual optimality and market clearing}
\label{sec:opt_mc}

Let $\bigl( \Omega, \fcal, (\fcal_t), \PR \bigr)$ be a 
filtered probability space supporting a standard
one-dimensional $(\fcal_t)$-Brownian motion $W$,
and let $(I, \ical, \iota)$ be a $\sigma$-finite measure
space.
We will use $W$ to model the exogenous demand for
a certain good in an economy and, ultimately, the
good's price process in equilibrium.
We use $I$ to index the producers in the sector of
the economy producing the given good.
We assume that all of the producers are {\em small\/}, in the
sense that their individual decisions do not impact prices.
In other words, we assume that all of the producers
are {\em price takers\/}.
In such a {\em competitive\/} setting, we assume that
the measure $\iota$ on $(I, \ical)$ is non-atomic.
In particular, we assume that $\{ i \} \in \ical$ and
$\iota \bigl( \{ i \} \bigr) = 0$ for all $i \in I$.

At time $0-$, the individual producer $i \in I$ has
installed capability $c_i > 0$.
We use the term ``capability'' to refer to a producer's
productive efficiency, which influences its production
costs as described below.

We denote by $P_t$ the price at which the
good is traded in the economy at time
$t \geq 0$.
In the competitive setting that we consider, each
individual producer takes this price as exogenously
given and has no individual influence over it.
Indeed, the price process $P$ is the {\em aggregate
outcome\/} of all producers' actions in equilibrium.
We also denote by $C_i (t) > 0$ and $Q_i (t)
\geq 0$ the capability and the rate of production
that the individual producer $i \in I$ maintains at
time $t>0$.
We assume that, over a small time interval
$[t, t {+} \Delta t [$, the production cost for
producer $i \in I$  is proportional to $C_i^{- a_i} (t)
Q_i^{1 + b} (t) \, \Delta t$, for some constants $b > 0$
and $a_i \in \mbox{} ]0, b[$.
To simplify the notation in our analysis, we re-parametrise 
and we assume that the cost of production incurred by
producer $i \in I$ during the time interval
$[t, t {+} \Delta t [$ is given by
\be
\frac{\beta}{\lambda _i^{1 / \beta} (1+\beta)
^{1 + 1/\beta}} \frac{Q_i^{1 + 1/\beta} (t)}
{C_i^{\alpha_i / \beta} (t)} \, \Delta t ,
\ee
for some constants $\beta > 0$, $\alpha_i \in
\mbox{} ]0, 1[$ and $\lambda_i > 0$.
In this setting, the production cost elasticity
$b = 1 / \beta$ is the {\em only\/} parameter
that is common among all heterogeneous
producers.

At any time $t \geq 0$, producer $i \in I$ determines
both the quantity $Q_i (t) \, \Delta t$ of the good
they produce and their expansion $\Delta C_i (t)$
over the time interval $[t, t {+} \Delta t[$.
Over the same time interval, the producer receives
the payoff
\be
\biggl( Q_i (t) P_t - \frac{\beta}
{\lambda_i^{1 / \beta} (1+\beta)^{1 + 1 / \beta}}
\frac{Q_i^{1 + 1/\beta} (t)}{C_i^{\alpha_i / \beta} (t)}
\biggr) \, \Delta t - \lambda_i k_i \, \Delta C_i (t) .
\ee
Here, the first term is the revenue from sales after
deducting production costs, while the constant
$\lambda_i k_i > 0$  in the second term is the cost
of unit expansion.
This expression implies that the optimal for producer
$i \in I$ production rate is
\ben
Q_i^\star (t) = \lambda_i (1+\beta) C_i^{\alpha_i}
(t) P_t^\beta , \label{Q_star}
\een
which gives rise to the payoff $C_i^{\alpha_i} (t)
P_t^{1+\beta} \Delta t$.
Accordingly, it is optimal for the individual producer
$i \in I$ to adopt an expansion strategy $C_i$
that maximises the performance index
\ben
\widetilde{J}_{c_i}^{(i)} (C_i \mid P) = \EXP \biggl[
\int _0^\infty e^{-r_i t} C_i^{\alpha_i} (t) P_t^{1+\beta}
\, \di t - k_i \int _{[0,\infty[} e^{-r_i t} \, \di C_i (t) \biggr] ,
\label{Ji-tilde}
\een
where $r_i > 0$ is the producer's discounting rate.

We assume that, over a small time interval
$[t, t {+} \Delta t [$, the demand for the good
is given by $D_t P_t^{-\delta} \, \Delta t$, where
$\delta > 0$ is the demand's price elasticity and the
process $D$ represents the exogenous base demand
corresponding to a unit price.
We model $D$ by
\be
D_t = D_0 \exp \bigl( \widetilde{\mu} t + \widetilde{\sigma}
W_t \bigr) ,
\ee
for some constants $D_0 > 0$, $\widetilde{\mu} \in \bbr$
and $\widetilde{\sigma} > 0$.
We determine the price process $P$ by imposing
the market-clearing condition, which requires that
demand should equal supply.
To this end, suppose that all producers have committed
on their individual expansion strategies $\{ C_i , \ i \in I\}$.
In view of (\ref{Q_star}), the market-clearing requirement
gives rise to the identity
\be
P_t^{-\delta} D_t = (1+\beta) P_t^\beta \int _I
\lambda _i  C_i^{\alpha_i} (t) \, \iota (\di i) ,
\ee
which implies that
\ben
P_t = X_t^{1/(1+\beta)} \Biggl( (1+\beta) \int _I
C_i^{\alpha_i} (t) \lambda _i \, \iota (\di i)
\Biggr) ^{-\gamma / (1+\beta)} , \label{Peq1}
\een
where $X = D^\gamma$, with
$\gamma = (1+\beta) / (\delta+\beta) \in
\mbox{} ]0,1[$.
In particular, $X$ is a geometric Brownian motion,
namely,
\ben
\di X_t = \mu X_t \, \di t + \sigma X_t \, \di W_t ,
\quad X_0 = D_0^\gamma , \label{dX}
\een
for some constants $\mu \in \bbr$ and $\sigma > 0$.

To ensure that the market-clearing condition
(\ref{Peq1}) is well-defined and associated with
a real-valued equilibrium price process $P$,
we make the following assumption.

\begin{ass} \label{A1} {\rm
The functions $I \ni i \mapsto c_i$, $I \ni i \mapsto \alpha_i$,
$I \ni i \mapsto \lambda_i$, $I \ni i \mapsto k_i$ and
$I \ni i \mapsto r_i$ are all strictly positive and
$\ical$-measurable.
Furthermore,
\ben
\kappa_0 := (1+\beta) \int _I \lambda _i c_i^{\alpha_i}
\, \iota (\di i) < \infty . \label{kappa0}
\een
} \end{ass}

\section{Competitive production equilibrium}
\label{sec:prodeq}

Before addressing the issue of equilibrium,
we note that the general solution to the ODE
\ben
\lscr_i u (x) := \frac{1}{2} \sigma^2 x^2 u'' (x)
+ \mu x u' (x) - r_i u (x) = 0 , \label{scriptL}
\een
which is associated with the geometric Brownian
motion $X$ defined by (\ref{dX}) and the discounting
rate $r_i > 0$ of producer~$i$, is given by
\ben
u(x) = \Delta_1 x^{m_i} + \Delta_2 x^{n_i},
\label{ODEhom-sol}
\een
for constants $\Delta_1 \in \bbr$ and $\Delta_2 \in \bbr$,
where $m_i < 0 < n_i$ are the solutions to the
quadratic equation
\be
\frac{1}{2} \sigma^2 \ell^2 + \left( \mu - \frac{1}{2}
\sigma ^2 \right) \ell - r_i = 0 ,
\ee
namely,
\be
m_i, n_i = \frac{- \bigl( \mu - \frac{1}{2} \sigma^2 \bigr) \mp
\sqrt{\left( \mu - \frac{1}{2} \sigma^2 \right) ^2 + 2 \sigma^2 r_i}}
{\sigma^2} .
\ee
For future reference, we note that the constants
$m_i < 0 < n_i$ are such that
\begin{gather}
r_i > \mu \ \Leftrightarrow \ n_i > 1 , \quad
n_i + m_i - 1 = - \frac{2 \mu}{\sigma^2} , \quad
n_i m_i = - \frac{2r_i}{\sigma^2} \label{rbnm1} \\
\text{and} \quad
r_i - \mu = \frac{1}{2} \sigma^2 (n_i - 1) (-m_i + 1) .
\label{rbnm2}
\end{gather}

We now make a further assumption about the data
of the optimisation problems (\ref{Ji-tilde}) faced
by the individual producers.
To this end, we define
\begin{equation}
\overline{\alpha} = \iota \textit{-} \esssup \alpha
= \inf \bigl\{ a \in \mbox{} ]0,1] \mid \
\iota (\alpha > a) = 0 \bigr\}.
\label{over-alpha} \\
\end{equation}

\begin{ass} \label{A2} {\rm
The functions $I \ni i \mapsto \alpha_i$ and
$I \ni i \mapsto n_i$ are such that%
\footnote{
Note that (\ref{rbnm1}) and the conditions in
(\ref{epsilon-i}) imply that
\ben
n_i > 1 \ \Leftrightarrow \ r_i > \mu \text{ for all }
i \in I . \label{ni>1}
\een
In the extremal case when $\overline{\alpha} = 1$,
the second condition in (\ref{epsilon-i}) simplifies to
(\ref{ni>1}).
}
\begin{equation}
\alpha_i \in \mbox{} ]0,1[ \quad \text{and} \quad
\varepsilon_i := \frac{(n_i - 1) (1-\alpha_i)}
{\alpha_i} - \frac{1-\overline{\alpha}}
{1-\overline{\alpha} + \overline{\alpha} \gamma}
> 0 \quad \text{for all } i \in I . \label{epsilon-i}
\end{equation}
Furthermore,
\be
\int _I \lambda_i \biggl( q
\frac{\alpha_i (n_i - 1)}{(r_i - \mu) n_i k_i} \biggr)
^{\alpha_i / (1-\alpha_i)} \, \iota (\di i) < \infty \quad
\text{for all constants } q > 0 .
\ee
} \end{ass}

\begin{de} \label{de:adm-contr} {\rm
The family $\acal$ of all admissible production
expansion strategies consists of all collections of
processes $\{ C_i , \ i \in I \}$ such that
\smallskip

\noindent (i)
the function $(\Omega \times \bbr_+) \times I \ni
\bigl( (\omega, t) , i \bigr) \mapsto C_i (\omega, t)$
is $\pcal \otimes \ical$-measurable, where $\pcal$
is the progressive $\sigma$-algebra on $\Omega
\times \bbr_+$, and
\smallskip

\noindent (ii)
given any $i \in I$, $C_i $ is an $(\fcal_t)$-adapted
increasing process with c\`{a}dl\`{a}g sample paths
such that $C_i (0-) = c_i$ and
\ben
\EXP \biggl[ \int _0^\infty e^{-r_i t} C_i^{\alpha_i} (t)
X_t \, \di t \biggr] < \infty . \label{Ac-IC}
\een

\noindent
Given $i \in I$, the family $\acal_i$ of all admissible 
{\em individual\/} expansion strategies consists of all
processes $C_i$ that are as in (ii).
} \end{de}

\noindent
The admissibility conditions (\ref{Ac-IC}) ensure
that, in equilibrium, the individual investors'
optimisation problems are well-posed.
Indeed, the market-clearing condition (\ref{Peq1})
and Assumption~\ref{A1} imply that
\be
P_t^{1+\beta} \leq \kappa_0^{-\gamma}
X_t \quad \text{for all } t \geq 0 .
\ee
Therefore, $\widetilde{J}_{c_i}^{(i)} (C_i \mid P)$
is well-defined for all $i \in I$.

\begin{de} \label{de:equil} {\rm
A competitive production equilibrium consists of a
continuous $(\fcal_t)$-adapted price process $P$
and a family of expansion strategies
$\{ C_i ^\star , \, i \in I \} \in \acal$ such that the
market-clearing condition (\ref{Peq1}) holds and
\ben
\widetilde{J}_{c_i}^{(i)} (C_i^\star \mid P) =
\sup _{C_i \in \acal_i} \widetilde{J}_{c_i}^{(i)}
(C_i \mid P) \quad \text{for all } i \in I ,
\label{Cistar-eq-de}
\een
where $\widetilde{J}_{c_i}^{(i)} (C_i \mid P)$
is defined by (\ref{Ji-tilde}).
} \end{de}

\section{The main result: a Markovian competitive
production equilibrium}
\label{sec:mr}

In the context of a Markovian equilibrium, intuition suggests
that  it should be optimal for each producer
to increase their capability only when the good's price
or demand is at a maximal level, namely, at times
$t \geq 0$ when $P_t = \overline{P}_t$ or $X_t =
\overline{X}_t$, where
\ben
\overline{P} _t = \sup _{0 \leq s \leq t} P_s
\quad \text{and} \quad
\overline{X} _t = \sup _{0 \leq s \leq t} X_s .
\label{PXbar}
\een
In light of this observation, we restrict attention
to expansion strategies of the form
\ben
C_i (t) = c_i \vee \Usym_i (\overline{P}_t)
\quad \text{and} \quad
C_i (t) = c_i \vee \gbG_i (\overline{X}_t)
, \label{C-PX}
\een
where $\bUsym = \{ \Usym_i, \ i \in I \}$ and
$\bgbG = \{ \gbG_i, \ i \in I \}$ are in the set $\bcal$.

\begin{de} \label{S-Bcal} {\rm
We define $\bcal$ to be the set of all families of
functions $\mathbb{X} = \{ \xcal_i, \ i \in I \}$ such
that\footnote{
Condition~(i) in this definition implies that the
collections $\{ \Usym_i (\overline{P}) , \ i \in I \}$
and $\{ \gbG_i (\overline{X}) , \ i \in I \}$ both
satisfy condition~(i) in Definition~\ref{de:adm-contr}.
}
\smallskip

\noindent (i)
the function $\bbr_+ \times I \ni (z,i) \mapsto
\mathbb{X} (z, i) = \xcal_i (z)$ is $\bcal (\bbr_+)
\otimes \ical$-measurable,
\smallskip

\noindent (ii)
the integrability condition
\ben
\int _I \lambda _i \xcal_i^{\alpha_i} (z)
\, \iota (\di i) < \infty \label{Chi-props2}
\een
holds true for all $z > 0$, and
\smallskip

\noindent (iii)
given any $i \in I$, $\xcal_i$ is the difference of two
convex functions\footnote{
A function $f: \bbr_+ \rightarrow \bbr$ is the difference
of two convex functions if and only if it is absolutely
continuous with left-hand derivative that is a function of
finite variation.
Given such a function, we denote by
$f_-'$ its left-hand side first derivative.
}
such that
\ben
\lim _{z \downarrow 0} \xcal_i (z) = 0 , \quad
\lim _{z \uparrow \infty} \xcal_i (z) = \infty
\quad \text{and} \quad
(\xcal_i) _-' (z) > 0 \text{ for all } z > 0 .
\label{Chi-props1}
\een
} \end{de}

In the context of (\ref{PXbar}) and (\ref{C-PX}),
expression (\ref{Peq1}) indicates that, under
market-clearing conditions, the good's equilibrium
price process $P$ should be characterised by the
expressions
\be
P^{1+\beta} = X \iscr [\bUsym] (\overline{P})
\quad \text{and} \quad
P^{1+\beta} = X \iscr [\bgbG] (\overline{X}) ,
\ee
where the operator $\iscr$ is defined by
\ben
\iscr [\mathbb{X}] (z) = \Biggl( (1+\beta) \int _I
\lambda _i \bigl( c_i^{\alpha_i} \vee \xcal_i^{\alpha_i}
(z) \bigr) \, \iota (\di i) \Biggr) ^{-\gamma} , \quad
z \geq 0 , \label{I[Xi]}
\een
for $\mathbb{X} = \{ \xcal_i, \ i \in I \} \in \bcal$.

Motivated by these observations, we
construct a competitive equilibrium starting from
a given family $\bUsym = \{ \Usym_i, \ i \in I \}
\in \bcal$.
With the possible exception of the producer indexed
by $i \in I$, we assume that every other producer,
say $j \in I \setminus \{ i \}$, adopts the capacity
expansion strategy
$C_j (t) = c_j \vee \Usym_j (\overline{P}_t)$.
In the competitive setting under consideration,
this specification gives rise to the price
process $P$ that solves the equation
\ben
P^{1+\beta} = X \iscr [\bUsym] (\overline{P})
=: X \usym (\overline{P}) . \label{eq.PbarP}
\een
Accordingly, the individual optimality requirement
(\ref{Cistar-eq-de}) reduces to determining an
expansion process $C_i^\star$ that maximises the
performance criterion
\ben
\EXP \biggl[ \int _0^\infty e^{-r_i t} C_i^{\alpha_i} (t)
\bigl( X_t \usym (\overline{P}_t) \bigr) ^{1+\beta}
\, \di t - k_i \int _{[0,\infty[} e^{-r_i t} \, \di C_i (t) \biggr]
\label{PI-psi}
\een
over all admissible $C_i$.

In Theorem~\ref{thm:Phi-Psi-eq}, we show that
the solution to equation (\ref{eq.PbarP}) can be
expressed as
\be
P^{1+\beta} = X \iscr [\bgbG] (\overline{X})
= X \hsym (\overline{X}) ,
\ee
where the family $\bgbG = \{ \gbG_i, \ i \in I \}$
and the function $\hsym$ are such that
\be
\gbG _i (\overline{x}) = \Usym _i \Bigl( \bigl( \overline{x}
\hsym (\overline{x}) \bigr) ^{1/(1+\beta)} \Bigr) .
\ee
Consequently, maximising (\ref{PI-psi}) is equivalent to
maximising the performance index
\ben
\EXP \biggl[ \int _0^\infty e^{-r_i t} C_i^{\alpha_i} (t)
\bigl( X_t \hsym (\overline{X}_t) \bigr) ^{1+\beta}
\, \di t - k_i \int _{[0,\infty[} e^{-r_i t} \, \di C_i (t) \biggr]
. \label{PI-phi}
\een

We next turn to the individual optimisation problem
faced by the singled-out producer $i \in I$.
Specifically, we solve the problem of maximising the
performance criterion (\ref{PI-phi}) over all admissible
expansion strategies, as detailed in
Sections~\ref{sec:singular-HJB}--\ref{sec:sol}.
The resulting solution shows that the producer's
optimal expansion strategy is given by
\ben
C_i^\dagger (t) = c_i \vee \Usym_i^\star \bigl(
\overline{P}_t \bigr) := c_i \vee \biggl(
\frac{\alpha_i (n_i - 1)}{(r_i - \mu) n_i k_i} \biggr)
^{1/(1-\alpha_i)} \overline{P}_t ^{(1+\beta) / (1-\alpha_i)}
.  \label{eq.C*0}
\een
Notably, this result {\em does not\/}
depend on the initial specification of $\usym =
\iscr [\bUsym]$ in (\ref{eq.PbarP}).
It follows that the solution to the equation
\ben
P^\star = \Bigl( X \iscr [\bUsym ^\star] (\overline{P}
^\star) \Bigr) ^{1 / (1+\beta)} , \label{eq.PbarP*}
\een
where $\bUsym ^\star = \{ \Usym_i^\star , \
i \in I \}$ are defined in (\ref{eq.C*0}), and the
production expansion strategies
\ben
C_i^\star (t) = c_i \vee \Usym_i^\star \bigl(
\overline{P}_t^\star \bigr) = c_i \vee \biggl(
\frac{\alpha_i (n_i - 1)}{(r_i - \mu) n_i k_i} \biggr)
^{1/(1-\alpha_i)} {\overline{P}_t^\star}
^{(1+\beta) / (1-\alpha_i)} \label{eq.C*}
\een
comprise a competitive production equilibrium
in {\em feedback\/} form.
Furthermore, this production equilibrium admits
the {\em closed\/} form representation given by
the expressions
\begin{gather}
P_t^\star = \Bigl( X_t \, \iscr [\bgbG^\star] (\overline{X}_t)
\Bigr) ^{1/(1+\beta)} \label{eq-barX1} \\
\text{and} \quad
C_i^\star (t) = c_i \vee \gbG_i^\star (\overline{X}_t)
:= c_i \vee \Usym_i^\star \bigl( \overline{X}_t
\, \hsym (\overline{X}_t) \bigr) ^{1 / (1-\alpha_i)} ,
\label{eq-barX2}
\end{gather}
where $\hsym$ is defined as in (\ref{h2phiPhii}).

We prove the next theorem, which is the main result
of the paper, in Section~\ref{sec:main-proof}.

\begin{thm} \label{thm:main}
Suppose that Assumptions~\ref{A1} and \ref{A2}
hold true.
A competitive production equilibrium in the sense of
Definition~\ref{de:equil} is given by
(\ref{eq.PbarP*}) and (\ref{eq.C*}), or equivalently,
by (\ref{eq-barX1}) and (\ref{eq-barX2}).
Furthermore, $\iscr [\bUsym^\star]$ and $\iscr
[\bgbG^\star]$, where $\iscr$ is defined by (\ref{I[Xi]}),
are differences of two convex functions
and the stochastic dynamics of the equilibrium price
process $P$ are given by
\begin{align}
\di \ln P_t^\star & = \frac{\mu - \frac{1}{2} \sigma^2}
{1+\beta} \, \di t + \frac{1}{1+\beta}
\frac{(\iscr [\bUsym^\star]) _-' (\overline{P}_t^\star)}
{\iscr [\bUsym^\star] (\overline{P}_t)} \, \di
\overline{P}_t^\star + \frac{\sigma}{1+\beta} \, \di W_t
\nonumber \\
& = \frac{\mu - \frac{1}{2} \sigma^2}{1+\beta} \,
\di t + \frac{1}{1+\beta}
\frac{(\iscr [\bgbG^\star]) _-' (\overline{X}_t)}
{\iscr [\bgbG^\star] (\overline{X}_t)} \, \di \overline{X}_t
+ \frac{\sigma}{1+\beta} \, \di W_t , \label{dP-eq}
\end{align}
with initial condition
\ben
P_0^\star = \Bigl( X_0 \, \iscr [\bgbG^\star] (X_0) \Bigr)
^{1/(1+\beta)} . \label{P0-eq}
\een
\end{thm}

\begin{rem} {\rm
The optimal expansion strategies $\bUsym ^\star$
given by (\ref{eq.C*}) are {\em independent\/} of the
original choice of $\bUsym$, in particular, of the
specific function $\usym = \iscr [\bUsym]$ used in
(\ref{eq.PbarP}).
This observation indicates that the form of competitive
equilibrium that we derive is rather {\em robust\/}.
Indeed, even if the prevailing price process $P$
differs from the equilibrium one $P^\star$,
the individual producers will still find it optimal to
adopt the expansion strategies given by
(\ref{eq.C*}), provided that $P$  conforms to the
structure given in (\ref{eq.PbarP}).
} \end{rem}

\section{The solution to equation (\ref{eq.PbarP})}
\label{sec:equiv}

\begin{lem} \label{I[Xi]-d2cf}
Consider any collection $\mathbb{X} = \{ \xcal_i,
\ i \in I \} \in \bcal$ and the operator $\iscr$ that is
defined by (\ref{I[Xi]}).
The function $\chi := \iscr [\mathbb{X}]$ is a difference
of two convex functions,
\begin{gather}
\lim _{z \downarrow 0} \chi (z) = \kappa_0^{-\gamma}
, \quad \lim _{z \uparrow \infty} \chi (z) = 0 ,
\label{I[Xi]-prop1} \\
\chi _-' (z) = - \gamma (1+\beta) \chi
^{(1+\gamma) / \gamma} (z)  \int _I \lambda_i
\alpha_i {\bf 1} _{\{ c_i < \xcal_i (z) \}} 
\xcal_i^{\alpha_i-1} (z) (\xcal_i) _-' (z) \, \iota
(\di i) \label{I[Xi]'} \\
\text{and} \quad 
\chi _-' (z) < 0 \text{ for all } z > \sup \big\{ z \geq 0
\mid \ \chi (z) = \kappa_0^{-\gamma} \bigr\} ,
\label{I[Xi]-prop2}
\end{gather}
where $\kappa_0$ is defined by (\ref{kappa0}).
\end{lem}
\noindent {\bf Proof.}
Each of the functions $\xcal_i^{\alpha_i}$ and
$c_i^{\alpha_i} \vee \xcal_i^{\alpha_i}$ is a
difference of two convex functions.
Using Fubini's theorem, we can therefore
see that
\begin{align}
\frac{1}{(1+\beta)} \bigl( \chi ^{-1/\gamma} (z_2)
- \chi ^{-1/\gamma}  (z_1) \bigr) & = \int _I
\int _{z_1}^{z_2} \lambda _i \bigl( c_i^{\alpha_i}
\vee \xcal_i^{\alpha_i} \bigr) _-' (y) \, \di y \, \iota
(\di i) \nonumber \\
& = \int _{z_1}^{z_2} \int _I \lambda _i \bigl(
c_i^{\alpha_i} \vee \xcal_i ^{\alpha_i} \bigr) _-' (y)
\, \iota (\di i) \, \di y \quad \text{for all } z_1 < z_2 ,
\nonumber
\end{align}
which implies that $\chi$ is absolutely continuous
with left-hand side derivative given by
\be
\chi _-' (z) = - \gamma (1+\beta)
\chi ^{(1+\gamma) / \gamma} (z) \int _I \lambda _i
\bigl( c_i^{\alpha_i} \vee \xcal_i^{\alpha_i}
\bigr) _-' (z) \, \iota (\di i) ,
\ee
namely, (\ref{I[Xi]'}).
The function $\chi _-'$ is of finite variation because this
is true for all of the functions $\bigl( c_i^{\alpha_i} \vee
\xcal_i^{\alpha_i} \bigr) _-'$.
Furthermore, (\ref{I[Xi]-prop1}) and (\ref{I[Xi]-prop2})
follow immediately from (\ref{I[Xi]}), (\ref{I[Xi]'})
and the fact that the functions $\xcal_i$ are strictly
increasing.
\mbox{}\hfill$\Box$

\begin{thm} \label{thm:Phi-Psi-eq}
Consider a family $\bUsym = \{ \Usym_i , \ i \in I \}
\in \bcal$ and let $X$ be a strictly positive
$(\fcal_t)$-progressive process.
Also, define
\ben
\hfn (\overline{p}) = \frac{\overline{p} ^{1+\beta}}
{\iscr [\bUsym] (\overline{p})} , \quad
\hsym (\overline{x}) = \bigl( \iscr [\bUsym] \circ
\hfn^\inv \bigr) (\overline{x})
\quad \text{and} \quad
\gbG _i (\overline{x}) = \Usym _i \Bigl( \bigl(
\overline{x} \hsym (\overline{x}) \bigr) ^{1/(1+\beta)}
\Bigr) , \label{h2phiPhii}
\een
where $h^\inv$ denoting the inverse of $\hfn$.
The following statements hold true.
\smallskip

\noindent
{\rm (I)}
The family $\bgbG = \{ \gbG_i , \ i \in I \}$ belongs
to $\bcal$.
Furthermore, the function $\hsym$ is such that
$\hsym = \iscr [\bgbG]$,
\ben
\bigl( \overline{x} \hsym (\overline{x}) \bigr)
_-' > 0 \quad \text{and} \quad
\lim _{\overline{x} \uparrow \infty} \overline{x}
\hsym (\overline{x}) = \infty . \label{phi-props}
\een

\noindent
{\rm (II)}
There exists a unique strictly positive
$(\fcal_t)$-progressive process $P$ such that
\ben
P^{1+\beta} = X \iscr [\bUsym] (\overline{P}) ,
\label{P-eqn}
\een
where $\overline{P} _t = \sup _{0 \leq s \leq t} P_s$.
This process is given by
\ben
P = \bigl( X \iscr [\bgbG] (\overline{X}) \bigr)
^{1/(1+\beta)} , \label{P-eqn-sol}
\een
where $\overline{X} _t = \sup _{0 \leq s \leq t} X_s$.
\end{thm}
\noindent {\bf Proof.}
We first note that
\be
\lim _{\overline{p} \downarrow 0} \hfn
(\overline{p}) = 0 , \quad
\hfn _-' (\overline{p}) = \bigl( \overline{p}
^{1+\beta} / \iscr [\bUsym] (\overline{p}) \bigr)
_-' > 1 \quad \text{and} \quad
\lim _{\overline{p} \uparrow \infty} \hfn
(\overline{p}) = \infty
\ee
because $\iscr [\bUsym]$ satisfies (\ref{I[Xi]-prop1})
and (\ref{I[Xi]-prop2}).
The claims in part~(I) of the theorem follow from
these properties of $\hfn$, the assumption that
$\bUsym \in \bcal$ and the equivalences
\begin{align}
\Bigl( \hfn (\overline{p}) \, (\iscr [\bUsym]
\circ \hfn ^\inv) \bigl( \hfn (\overline{p}) \bigr)
\Bigr) ^{1 / (1+\beta)} = \overline{p}
\quad & \Leftrightarrow \quad
\bigl( \overline{x} \, \bigl( \iscr [\bUsym] \circ
\hfn ^\inv \bigr) (\overline{x}) \bigr) ^{1 / (1+\beta)}
= \hfn ^\inv (\overline{x}) \nonumber \\
& \Leftrightarrow \quad
\hfn_1 (\overline{x}) :=
\bigl( \overline{x} \hsym (\overline{x}) \bigr)
^{1 / (1+\beta)} = \hfn ^\inv (\overline{x}) .
\nonumber
\end{align}
Furthermore, the third identity in these equivalences
and the expression
\be
\iscr [\bgbG] (\overline{x}) = \Biggl( (1+\beta)
\int _I \lambda _i \Bigl( c_i^{\alpha_i} \vee
\Usym _i^{\alpha_i} \Bigl( \bigl( \overline{x}
\hsym (\overline{x}) \bigr) ^{1 / (1+\beta)} \Bigr)
\Bigr) \, \iota (\di i) \Biggr) ^{-\gamma}
= \bigl( \iscr [\bUsym] \circ \hfn^\inv \bigr)
(\overline{x})
\ee
imply that
\ben
\iscr [\bUsym] \circ \hfn_1
= \iscr [\bUsym] \circ \hfn^\inv = \iscr [\bgbG]
= \hsym . \label{IPsi-h1h-IPhi}
\een

To establish part~(II) of the theorem, let $P$
be the process defined by (\ref{P-eqn-sol}).
This process is such that $\overline{P} = h_1
(\overline{X})$ because $\hfn_1$ is strictly
increasing (see (\ref{phi-props})).
However, this identity and (\ref{IPsi-h1h-IPhi})
imply that
\be
\iscr [\bgbG] (\overline{X}) =
\bigl( \iscr [\bUsym] \circ h_1 \bigr) (\overline{X}) =
\iscr [\bUsym] (\overline{P}) ,
\ee
which shows that $P$ is a solution to equation
(\ref{P-eqn}).

To prove uniqueness, suppose that $P$ is a solution
to equation (\ref{P-eqn}).
Such a process satisfies
\be
\overline{P}^{1+\beta} = \overline{X} \iscr [\bUsym]
(\overline{P})
\quad \Leftrightarrow \quad
\overline{P} = h^\inv (\overline{X}) .
\ee
However, these identities and (\ref{IPsi-h1h-IPhi})
imply that $\iscr [\bUsym] (\overline{P}) =
\iscr [\bgbG] (\overline{X})$.
It follows that every solution $P$ to equation
(\ref{P-eqn}) of the form given by
(\ref{P-eqn-sol}).
\mbox{}\hfill$\Box$

\section{Asymptotic results}
\label{sec:asymptotics}

\begin{lem} \label{lem:Kostas}
For any $\ical$-measurable functions $\eta : I
\rightarrow [0, \infty[$, $\zeta : I \rightarrow \mbox{}
]0, \infty[$ and $\xi : I \rightarrow \bbr$ such that
\be
\int _I \zeta_i \, \iota (\di i) < \infty
\quad \text{and} \quad
O(y) := \int _I {\bf 1} _{\{ \eta_i < y \}} \zeta_i
y^{\xi_i} \, \iota (\di i) < \infty \quad
\text{for all } y > 0,
\ee
it holds that
\ben
\lim _{y \uparrow \infty} \frac{\ln O(y)}{\ln y}
= \iota \textit{-} \esssup \xi =: \overline{\xi} .
\label{lnO-lim}
\een
\end{lem}
\noindent {\bf Proof.}
We first note that
\begin{align}
\limsup _{y \uparrow \infty} \frac{\ln O(y)}{\ln y}
& \leq \lim _{y \uparrow \infty} \frac{1}{\ln y} \ln
\int _I \zeta_i e^{\xi_i \ln y} \, \iota (\di i) =
\lim _{y \uparrow \infty} \frac{1}{\ln y} \ln
\int _I e^{\xi_i \ln y} \, \widetilde{\iota} (\di i)
\nonumber \\
& = \lim _{y \uparrow \infty} \ln {\bigl\lVert}
e^\xi \bigr\Vert _{L^{\ln y} (\widetilde{\iota})}
= \widetilde{\iota} \textit{-} \esssup \xi
= \iota \textit{-} \esssup \xi , \label{eq:Kostas1}
\end{align}
where $\widetilde{\iota}$ is the \emph{finite} measure on $(I, \ical)$
that is equivalent to $\iota$, with Radon-Nikodym
derivative given by $\di \widetilde{\iota} / \di \iota = \zeta$.

Next, for $k \in \mathbb{N}$ define
$I_k := \{ \eta < k \} \in \ical$, and note that
\begin{align}
\liminf _{y \uparrow \infty} \frac{\ln O(y)}{\ln y}
& \geq \lim _{y \uparrow \infty} \frac{1}{\ln y} \ln
\int _I {\bf 1}_{I_k} e^{\xi_i \ln y} \,
\widetilde{\iota} (\di i) \nonumber \\
& = \lim _{y \uparrow \infty} \ln {\bigl\lVert}
{\bf 1}_{I_k} e^\xi \bigr\Vert
_{L^{\ln y} (\widetilde{\iota})}
= \iota \textit{-} \esssup {\bigl( {\bf 1} _{I_k}
\xi \bigr)} . \nonumber
\end{align}
Combining this result with the fact that $\bigcup
_{k=1}^\infty I_k = I$, we obtain
\be
\liminf _{y \uparrow \infty} \frac{\ln O(y)}{\ln y}
\geq \iota \textit{-} \esssup \xi .
\ee
This last inequality and (\ref{eq:Kostas1})
imply (\ref{lnO-lim}).
\mbox{}\hfill$\Box$
\bigskip

Using this lemma and the equivalence
\ben
\lim _{y \uparrow \infty} \frac{\ln Q(y)}{\ln y}
= - \ell \quad \Leftrightarrow \quad
\lim _{y \uparrow \infty} y^\xi Q(y)
= \begin{cases} 0 , & \text{if } \xi < \ell , \\
\infty , & \text{if } \xi > \ell , \end{cases}
\label{lnQ/lnp-lims}
\een
where $Q$ is a strictly positive function,
we can establish the following asymptotic
result.

\begin{lem} \label{lem:est}
Suppose that Assumption~\ref{A1} holds true and
let $\bUsym = \{ \Usym_i , \ i \in I \} \in \bcal$ be
such that
\be
\Usym _i (\overline{p}) = K_i \, \overline{p}
^{(1+\beta) / (1-\alpha_i)} ,
\ee
for some constant $\beta > 0$ and some
$\ical$-measurable functions $K > 0$ and $\alpha$
such that $\alpha_i \in \mbox{} ]0,1[$ for all
$i \in I$.
Also, let $\hfn$, $\hsym$,  $\bgbG = \{ \gbG_i
, \ i \in I \}$ be as in Theorem~\ref{thm:Phi-Psi-eq}
and define $\usym = \iscr [\bUsym]$.
The functions $\hsym$ and $\usym$ are such that
\begin{gather}
\lim _{\overline{p} \uparrow \infty}
\frac{\ln \usym (\overline{p})}{\ln \overline{p}}
= - \frac{\overline{\alpha} (1+\beta) \gamma}
{1-\overline{\alpha}} , \quad
\lim _{\overline{p} \uparrow \infty}
\frac{\ln {\bigl( - \usym _-' (\overline{p}) \bigr)}}
{\ln \overline{p}} = -
\frac{\overline{\alpha} (1+\beta) \gamma}
{1-\overline{\alpha}} - 1 ,
\label{lnpsi-lims} \\
\lim _{\overline{x} \uparrow \infty}
\frac{\ln \hsym (\overline{x})}{\ln \overline{x}}
= - \frac{\overline{\alpha} \gamma}
{1-\overline{\alpha} + \overline{\alpha} \gamma} ,
\quad \lim _{\overline{x} \uparrow \infty}
\frac{\ln {\bigl( - \hsym _-' (\overline{x}) \bigr)}}
{\ln \overline{x}} =
- \frac{\overline{\alpha} \gamma}
{1-\overline{\alpha} + \overline{\alpha} \gamma}
- 1 , \label{lnphi-lims} \\
\lim _{\overline{x} \uparrow \infty}
\frac{\ln {\bigl(} \overline{x} \hsym (\overline{x}) \bigr)}
{\ln \overline{x}} = \frac{1-\overline{\alpha}}
{1-\overline{\alpha} + \overline{\alpha} \gamma}
\quad \text{and} \quad
\lim _{\overline{x} \uparrow \infty}
\frac{\ln {\bigl( x \hsym (\overline{x}) \bigr) _-'}}
{\ln \overline{x}} =
- \frac{\overline{\alpha} \gamma}
{1-\overline{\alpha} + \overline{\alpha} \gamma}
. \label{ln(xphi)-lims}
\end{gather}
where $\overline{\alpha}$ is defined by
(\ref{over-alpha}).
\end{lem}
\noindent {\bf Proof.}
In view of the definition of $\iscr [\bUsym]$,
we can see that
\be
(1+\beta) \int _I \lambda _i \Usym _i^{\alpha_i}
(\overline{p}) \, \iota (\di i) \leq \usym^{-1/\gamma}
(\overline{p}) \leq \kappa_0 + (1+\beta)
\int _I \lambda _i \Usym _i^{\alpha_i} (\overline{p})
\, \iota (\di i) ,
\ee
where $\kappa_0$ is defined by (\ref{kappa0}).
These inequalities, Lemma~\ref{lem:Kostas}
and the fact that the function $]0, 1[
\mbox{} \ni a \mapsto a / (1 - a)$
is increasing imply the first limit in
(\ref{lnpsi-lims}).
On the other hand, combining the expression
\begin{align}
\ln {\bigl( - \usym_-' (\overline{p}) \bigr)} =
\mbox{} & \ln \bigl( \gamma (1+\beta) \bigr) +
\frac{\gamma + 1}{\gamma} \ln \usym (\overline{p})
\nonumber \\
& + \ln \int _I {\bf 1}
_{\bigl\{ (c_i/K_i) ^{(1-\alpha_i) / (1+\beta)}
< \, \overline{p} \bigr\}}
\frac{\alpha_i (1+\beta) \lambda _i K_i^{\alpha_i}}
{1-\alpha_i}
\, \overline{p} ^{\alpha_i (1+\beta) / (1-\alpha_i) - 1}
\, \iota (\di i) , \nonumber
\end{align}
which follows from (\ref{I[Xi]'}), with
Lemma~\ref{lem:Kostas} and the first limit in
(\ref{lnpsi-lims}), we obtain the second limit in
(\ref{lnpsi-lims}).

Using the definitions of $\hfn$ and $\hsym$ as
well as (\ref{lnpsi-lims}), we can see
that
\begin{gather}
\lim _{\overline{p} \uparrow \infty}
\frac{\ln \hfn (\overline{p})}{\ln \overline{p}}
= 1+\beta - \lim _{\overline{p} \uparrow \infty}
\frac{\ln \usym (\overline{p})}{\ln \overline{p}}
= (1+\beta) \biggl( 1 +
\frac{\overline{\alpha} \gamma}{1-\overline{\alpha}}
\biggr) , \nonumber \\
\lim _{\overline{x} \uparrow \infty}
\frac{\ln \hfn^\inv (\overline{x})}{\ln \overline{x}}
= \lim _{\overline{x} \uparrow \infty}
\frac{\ln \hfn^\inv (\overline{x})}
{\ln \hfn \bigl( \hfn^\inv (\overline{x}) \bigr)}
= (1+\beta) ^{-1} \biggl( 1 +
\frac{\overline{\alpha} \gamma}{1-\overline{\alpha}}
\biggr) ^{-1} \nonumber \\
\text{and} \quad
\lim _{\overline{x} \uparrow \infty}
\frac{\ln \hsym (\overline{x})}{\ln \overline{x}}
= \lim _{\overline{x} \uparrow \infty}
\frac{\ln \usym \bigl( \hfn^\inv (\overline{x}) \bigr)}
{\ln \hfn^\inv (\overline{x})}
\frac{\ln \hfn^\inv (\overline{x})}{\ln \overline{x}}
= - \frac{\overline{\alpha} \gamma}
{1-\overline{\alpha} + \overline{\alpha} \gamma} .
\nonumber
\end{gather}
Furthermore,
\be
\lim _{\overline{p} \uparrow \infty}
\frac{\ln {\bigl( (\beta+1) \overline{p}^\beta
\usym ^{-1} (\overline{p}) \bigr)}}
{\ln \overline{p}} =
\lim _{\overline{p} \uparrow \infty}
\frac{\ln {\bigl( \overline{p}^{\beta+1} (-\usym)_-'
(\overline{p}) \usym^{-2} (\overline{p}) \bigr)}}
{\ln \overline{p}} = \beta +
\frac{\overline{\alpha} (1+\beta) \gamma}
{1-\overline{\alpha}} .
\ee
The last two of these limits and the equivalence
in (\ref{lnQ/lnp-lims}) imply that
\begin{align}
\lim _{\overline{p} \uparrow \infty}
\overline{p}^\xi \hfn _-' (\overline{p}) & =
\lim _{\overline{p} \uparrow \infty}
\overline{p}^\xi \biggl(
\frac{(\beta+1) \overline{p}^\beta}{\usym (\overline{p})}
+ \frac{\overline{p}^{\beta+1} (-\usym)_-' (\overline{p})}
{\usym ^2 (\overline{p})} \biggr) \nonumber \\
& = \begin{cases} 0 , & \text{if } \xi <
- (1+\beta) \bigl( 1 +
\frac{\overline{\alpha} \gamma}{1-\overline{\alpha}}
\bigr) + 1 , \\ \infty , & \text{if } \xi >
- (1+\beta) \bigl( 1 +
\frac{\overline{\alpha} \gamma}{1-\overline{\alpha}}
\bigr) + 1 . \end{cases} \nonumber
\end{align}
Combining this observation with (\ref{lnQ/lnp-lims}),
we obtain
\be
\lim _{\overline{p} \uparrow \infty}
\frac{\ln \hfn_-' (\overline{p})}{\ln \overline{p}}
= (1+\beta) \biggl( 1 +
\frac{\overline{\alpha} \gamma}{1-\overline{\alpha}}
\biggr) - 1 .
\ee
In view of the limits that we have derived thus
far and the calculation
\be
\lim _{\overline{x} \uparrow \infty}
\frac{\ln {\bigl( - \hsym _-' (\overline{x}) \bigr)}}
{\ln \overline{x}}
= \lim _{\overline{x} \uparrow \infty} \Biggl(
\frac{\ln {(- \usym)} _-' \bigl( \hfn^\inv (\overline{x}) \bigr)}
{\ln \hfn^\inv (\overline{x})}
- \frac{\ln \hfn_-' \bigl( \hfn^\inv (\overline{x}) \bigr)}
{\ln \hfn^\inv (\overline{x})} \Biggr)
\frac{\ln \hfn^\inv (\overline{x})}{\ln \overline{x}} ,
\ee
we can see that the second limit in (\ref{lnphi-lims})
also holds true.

Finally, the limits in (\ref{ln(xphi)-lims}) can be seen
using similar arguments.
\mbox{}\hfill$\Box$

\section{A singular stochastic control problem
and its HJB equation}
\label{sec:singular-HJB}

In view of  the analysis in Sections~\ref{sec:mr}
and~\ref{sec:equiv}, we now consider the problem of
maximising the performance index defined by
(\ref{PI-phi}) over all admissible $C_i$.
In particular, we determine $C_i^\dagger \in \acal_i$
that maximises the performance index
\ben
J_{c_i,x,\overline{x}}^{(i)} (C) = \EXP \biggl[
\int _0^\infty e^{-r_i t} C^{\alpha_i} (t) X_t
\hsym (\widetilde{X}_t) \, \di t - k_i \int _{[0,\infty[}
e^{-r_i t} \, \di C(t) \biggr] , \label{Ji}
\een
for $x = \overline{x}$, where $X$ is given by (\ref{dX})
and
\ben
\widetilde{X}_t = \overline{x} \vee \overline{X}_t
= \overline{x} \vee \sup _{0 \leq s \leq t} X_s .
\label{Xtilde}
\een
Here, we allow for any initial condition
$\overline{x} \geq x$ because we will solve the
problem using dynamic programming.
Additionally, $\hsym$ is {\em any\/} function satisfying
the conditions of following assumption.

\begin{ass} \label{A4} {\rm
The function $\hsym$ can be expressed as a difference
of two convex functions and has the properties listed
in (\ref{phi-props}).
Moreover, for every $\varepsilon > 0$, there
exist constants
$\widetilde{K} (\varepsilon)$ and
$\widetilde{x} (\varepsilon) > 1$
such that
\ben
\overline{x} \hsym (\overline{x}) \leq \widetilde{K}
(\varepsilon) \, \overline{x}
^{( 1 - \overline{\alpha}) / (1-\overline{\alpha}
+ \overline{\alpha} \gamma) + \varepsilon}
\quad \text{for all } \overline{x} > \widetilde{x}
(\varepsilon) . \label{phi-asym-ass}
\een
} \end{ass}

We are thus faced with the singular stochastic control
problem whose value function is defined by
\ben
v^{(i)} (c_i, x, \overline{x}) = \sup _{C \in \acal_i}
J_{c_i,x,\overline{x}}^{(i)} (C) , \label{IP-v}
\een
where the class of admissible controls $\acal_i$
is specified by Definition~\ref{de:adm-contr}.
The HJB equation of this control problem is given by
\ben
\max \bigl\{ \lscr_i w (c, x, \overline{x}) + c^{\alpha_i} x
\hsym (\overline{x}) , \ w_c (c, x, \overline{x}) - k_i \bigr\}
= 0 , \label{HJB}
\een
with boundary condition
\ben
w_{\overline{x}} (c, x, \overline{x}) = 0 , \label{HJB-BC}
\een
where the operator $\lscr_i$ is defined by (\ref{scriptL}).

We will show that the the value function $v^{(i)}$
identifies with the solution $w = w^{(i)}$ to the
HJB equation (\ref{HJB}) with boundary condition
(\ref{HJB-BC}) that is of the form presented by
Problem~\ref{FBP} below.
This involves a function $G_i : \bbr_+^2
\rightarrow \bbr_+$ such that $G_i (\cdot,
\overline{x})$ is $C^1$,
\ben
(G_i) _c (c,\overline{x}) > 0 , \quad
\lim _{c \downarrow 0} G_i (c,\overline{x}) = 0
\quad \text{and} \quad
\lim _{c \uparrow 0} G_i (c,\overline{x}) =  \infty
\quad \text{for all } \overline{x} > 0 . \label{G-prop}
\een
It also involves the unique solutions $\gbG_i
(\overline{x})$ to the equations
\ben
G_i \bigl( \gbG_i (\overline{x}) , \overline{x} \bigr) =
\overline{x} , \quad \overline{x} > 0 , \label{gG-G-eqn}
\een
which give rise to a function $\gbG_i$ that is required
to be the difference of two convex functions and
such that
\ben
(\gbG_i) _-' (\overline{x}) > 0 \text{ for all } \overline{x}
> 0 , \quad
\lim _{\overline{x} \downarrow 0} \gbG_i (\overline{x})
= 0 \quad \text{and} \quad
\lim _{\overline{x} \uparrow \infty} \gbG_i (\overline{x})
= \infty . \label{gbG-lims}
\een

\begin{figure}[!t]
  \centering
  \begin{minipage}[b]{0.49\textwidth}
    \includegraphics[width=\textwidth]{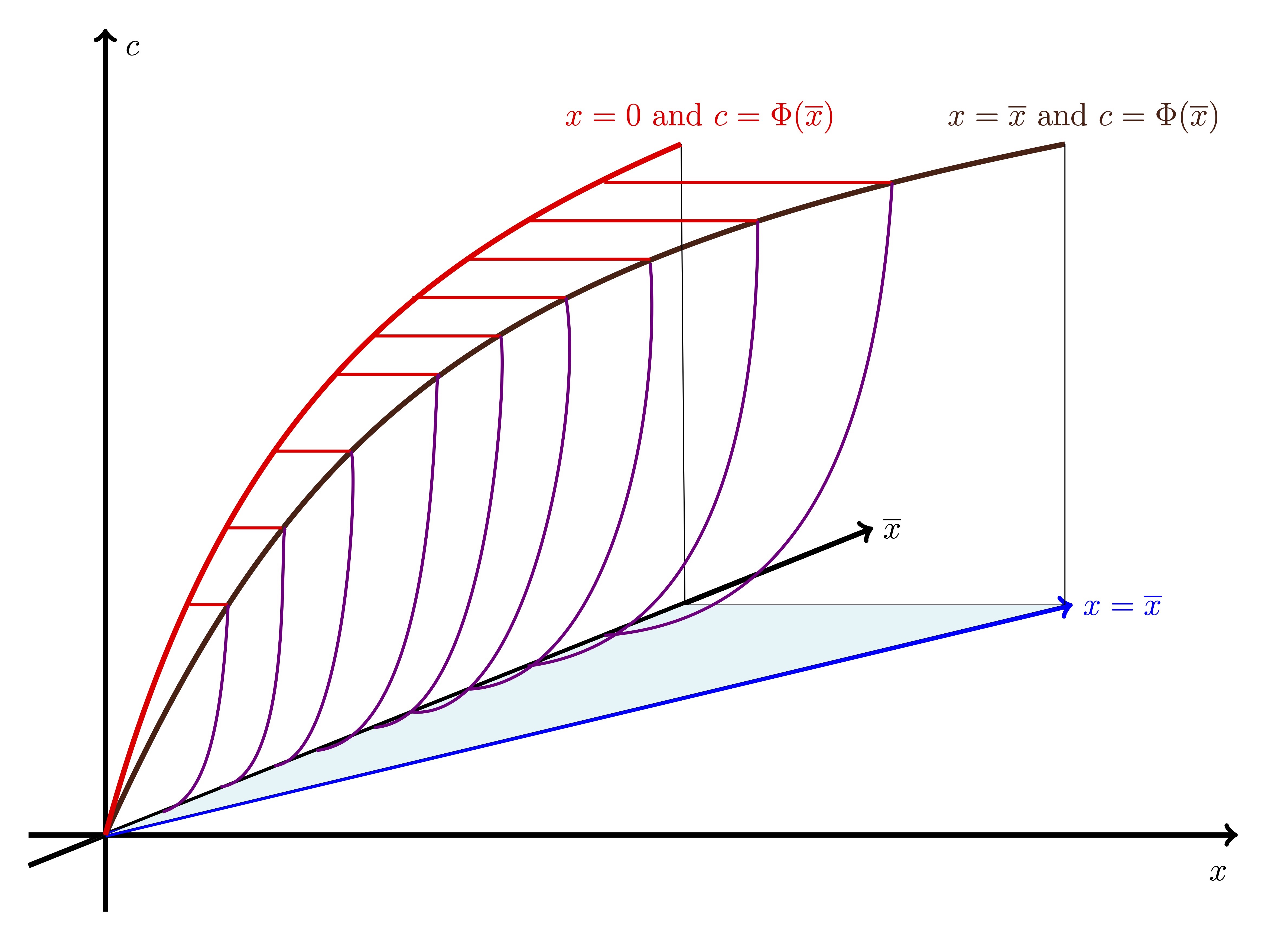}
    \caption{\small Graph illustrating the surfaces
    $\scal_i^1$ (purple lines) and $\scal_i^2$ (red lines).}
    \label{pic1}
  \end{minipage}
  \hfill
  \begin{minipage}[b]{0.49\textwidth}
    \includegraphics[width=\textwidth]{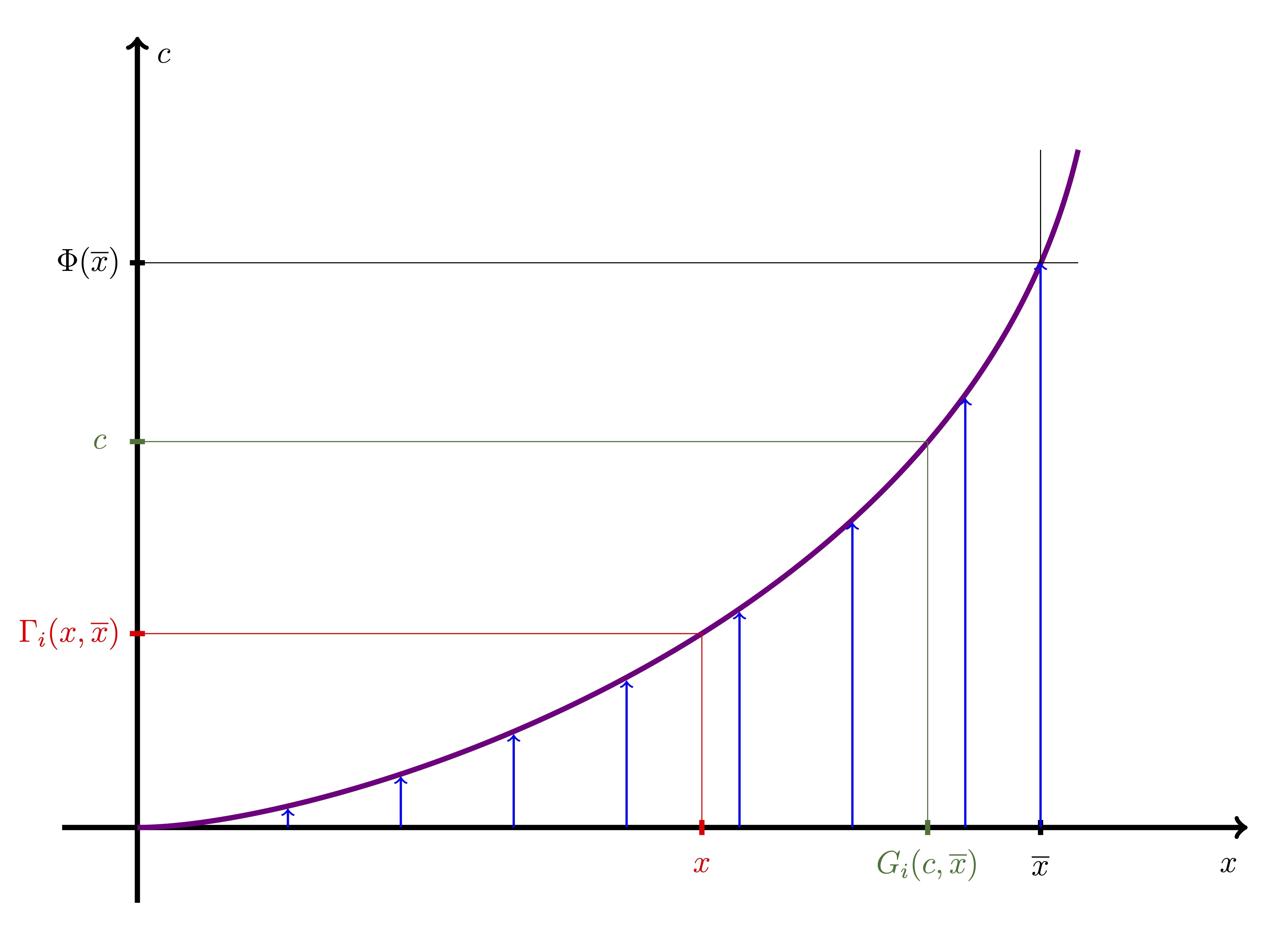}
    \caption{\small Graph illustrating the functions $G_i$ and
    $\Gamma_i$.}
    \label{pic2}
  \end{minipage}
\end{figure}

The functions $G$ and $\gbG$ define the free-boundary
surfaces
\begin{align}
\scal_i^1 & = \bigl\{ (c, x, \overline{x}) \in \bbr_+ \mid \
0 < x \leq \overline{x} \text{ and } x = G_i (c, \overline{x})
\bigr\} \label{surf1} \\
\text{and} \quad
\scal_i^2 & = \bigl\{ (c, x, \overline{x}) \in \bbr_+ \mid \
0 < x \leq \overline{x} \text{ and } c = \gbG_i
(\overline{x}) \bigr\} \label{surf2}
\end{align}
that can be illustrated by Figures~\ref{pic1}
and~\ref{pic2}.
These surfaces partition the control problem's state
space
\ben
\ccal = \bigl\{ (c, x, \overline{x}) \in \bbr_+^3 \mid
\ 0 < x \leq \overline{x} \bigr\} \label{C-st.sp.}
\een
into the sets
\begin{align}
\ccal _i^1 & = \bigl\{ (c, x, \overline{x}) \in \bbr_+ \mid \
0 < x \leq  \overline{x} \text{ and } c \leq \Gamma_i
(x, \overline{x}) \bigr\} , \label{C1} \\
\ccal _i^2 & = \bigl\{ (c, x, \overline{x}) \in \bbr_+ \mid \
0 < x \leq \overline{x} \text{ and } \Gamma _i
(x, \overline{x}) < c < \gbG _i (\overline{x}) \bigr\}
\label{C2} \\
\text{and} \quad
\ccal _i^3 & = \bigl\{ (c, x, \overline{x}) \in \bbr_+ \mid \
0 < x \leq \overline{x} \text{ and } \gbG _i (\overline{x})
\leq c \bigr\} , \label{C3} \\
\end{align}
where
\ben
\Gamma_i (\cdot, \overline{x}) = G_i^{\inv} (\cdot,
\overline{x}) , \label{Gamma}
\een
in which expression, $G_i^{\inv} (\cdot, \overline{x})$
the inverse of the function $G_i^{\inv} (\cdot,
\overline{x})$ (see also Figure~\ref{pic2}).

We postulate that the optimal expansion
strategy can be described informally as
follows.
At time 0, if the initial state $(c, x, \overline{x})$
is originally inside the region $\ccal _i^1$, then it is
optimal to exert minimal action so that the state process
is repositioned on the boundary surface
$\scal_i^2$.
In view of standard singular stochastic control
theory, such an action is associated with the
requirement that candidate $w$ for the value
function $v^{(i)}$ should satisfy (\ref{vc-eq})
as well as the inequality in (\ref{v-ODE}).
Beyond such a possible jump at time 0, it is optimal
to exert minimal action so that the state process
is prevented from entering the interior
of $\ccal _i^1$.
On the other hand, it is optimal to exert no control
effort while the state process takes values in the
interior of the set $\ccal _i^2 \cup \ccal _i^3$,
which is associated with the inequalities
(\ref{vc-ineq-u}) and (\ref{vc-ineq-d}),
as well as the equalities in (\ref{v-ODE}).
The significance of the surface $\scal_i^2$ arises
from the fact that, eventually, it is optimal to exert
minimal effort so as to prevent the state process
falling below the curve defined by $x = \overline{x}$
and $c = \gbG_i (\overline{x})$, which is the
intersection of $\scal_i^2$ with the boundary
of the state space defined by $x = \overline{x}$.
In particular, it is optimal that minimal control effort
should be exercised so that the state process
takes values on the surface $\scal_i^2$ at all times
after this surface has been reached.

\begin{problem} \label{FBP}
Determine a function $G_i: \bbr_+^2 \rightarrow \bbr_+$
and a function $w = w^{(i)} : \ccal \rightarrow \bbr$,
satisfying the following conditions.
\smallskip

\noindent {\rm (I)}
The function $G_i (\cdot, \overline{x})$ is $C^1$,
while the functions $G_i (c, \cdot)$ and $\gbG_i$
are differences of two convex functions.
Furthermore, $G_i$ and $\gbG_i$ satisfy
(\ref{G-prop})--(\ref{gbG-lims}).
\smallskip

\noindent {\rm (II)}
The function $w (\cdot, \cdot, \overline{x})$ is
$C^{1,2}$ in the interior of $\ccal$ for all
$\overline{x} > 0$.
\smallskip

\noindent {\rm (III)}
Given any $c \geq 0$ and $x>0$, the function
$w \bigl( c, x, \cdot \bigr) : \mbox{}
]x, \infty[ \mbox{} \rightarrow \bbr$ is the
difference of two convex functions.
Furthermore,
\ben
w_{\overline{x}} \bigl( c, \overline{x}, \overline{x} \bigr)
:= \lim _{x \uparrow \overline{x}} w_{\overline{x}-}
\bigl( c, x, \overline{x} \bigr) = 0 \quad \text{for all }
0 < \gbG_i (\overline{x}) \leq c . \label{vz-eq}
\een

\noindent {\rm (IV)}
The function $w$ is such that
\begin{align}
w_c (c, x, \overline{x}) & < k_i \quad \text{for all }
0 < x < \overline{x} \text{ and } \gbG_i (\overline{x})
\leq c , \label{vc-ineq-u} \\
w_c (c, x, \overline{x}) & < k_i \quad \text{for all }
0 < x \leq \overline{x} \text{ and } \Gamma _i
(x, \overline{x}) < c < \gbG _i (\overline{x}) ,
\label{vc-ineq-d} \\
w_c (c, x, \overline{x}) & = k_i \quad \text{for all }
0 < x \leq \overline{x} \text{ and } c \leq \Gamma_i
(x, \overline{x}) , \label{vc-eq}
\end{align}
and
\begin{align}
\lscr_i w (c, x, \overline{x}) + c^{\alpha_i} x
\hsym (\overline{x})
\begin{cases}
= 0 & \text{for all } 0 < x < \overline{x} \text{ and }
0 < \gbG_i (\overline{x}) \leq c , \\
= 0 & \text{for all } 0 < x \leq \overline{x} \text{ and }
\Gamma _i (x, \overline{x}) < c < \gbG _i
(\overline{x}) , \\
< 0 & \text{for all } 0 < x < \overline{x} \text{ and }
c \leq \Gamma_i (x, \overline{x}) ,
\end{cases} \label{v-ODE}
\end{align}
where $\lscr_i$ is defined by (\ref{scriptL}).
\end{problem}

\section{The solution to the control problem's
HJB equation}
\label{sec:HJB-sol}

In this section, we construct a solution $w =
w^{(i)}$ to the HJB equation (\ref{HJB}) with
boundary condition (\ref{HJB-BC}) that is as in
Problem~\ref{FBP}
and identifies with the value function $v^{(i)}$
defined by (\ref{IP-v}).
To this end, we consider the function
\begin{align}
w(c, x, \overline{x}) = \begin{cases}
c^{\alpha_i} \bigl( A_i (\overline{x}) + \int
_{\gbG_i (\overline{x})}^c f_i (y) \, \di y \bigr) x^{n_i}
+ \frac{1}{r_i - \mu} c^{\alpha_i} x \hsym (\overline{x}) ,
& \text{if } (c, x, \overline{x}) \in \ccal _i^3 , \\
B_i (c, \overline{x}) x^{n_i} + \frac{1}{r_i - \mu} c^{\alpha_i}
x \hsym (\overline{x}) , & \text{if } (c, x, \overline{x})
\in \ccal _i^2 \cup \scal _i^1 , \\
w \bigl( \Gamma_i (x, \overline{x}) , x, \overline{x}
\bigr) - k_i \bigl( \Gamma_i (x, \overline{x}) - c \bigr)
, & \text{if } (c, x, \overline{x}) \in \ccal _i^1
\setminus \scal _i^1 , \end{cases}
\label{v-expr}
\end{align}
where $\scal _i^2$, $\ccal _i^1$, $\ccal _i^2$,
$\ccal _i^3$ and $\Gamma_i$ are given by
(\ref{surf2}), (\ref{C1})--(\ref{C3}) and (\ref{Gamma}).
In view of the general solution (\ref{ODEhom-sol})
to the ODE $\lscr_i u(x) = 0$, we have taken the
coefficient of $x^{m_i}$ to be 0 here because,
otherwise, the function $w$ would not satisfy the
transversality condition that is required for its
identification with the control problem's value
function.

To determine the function $B_i$ and the free-boundary
function $G_i$, we appeal to the so called
``smooth-pasting condition'' of singular stochastic
control.
In particular,  we require that $w (\cdot, x, \overline{x})$
should be $C^2$ along the free-boundary point
$G_i (c,\overline{x})$, which suggests the equations
\begin{align}
\lim _{x \uparrow G_i  (c, \overline{x})} w_c
(c,x,\overline{x}) = (B_i)_c (c, \overline{x}) G_i^{n_i}
(c,\overline{x}) + \frac{\alpha_i}{r_i - \mu} c^{\alpha_i - 1}
\hsym (\overline{x}) G_i  (c, \overline{x}) & = k_i
\label{v_c-G} \\
\text{and} \quad
\lim _{x \uparrow G_i (c, \overline{x})} w_{cx}
(c,x,\overline{x}) = n_i (B_i)_c (c, \overline{x}) G_i^{n_i - 1}
(c,\overline{x}) + \frac{\alpha_i}{r_i - \mu} c^{\alpha_i - 1}
\hsym (\overline{x}) & = 0 . \label{v_cx-G}
\end{align}
The solution to this system of equations is given by
\ben
G_i (c,\overline{x}) = \frac{(r_i - \mu) n_i k_i}{\alpha_i (n_i - 1)}
\frac{c^{1-\alpha_i}}{\hsym (\overline{x})}
\quad \text{and} \quad
(B_i)_c (c,\overline{x}) = - \frac{k_i}{n_i - 1} G_i^{-n_i}
(c,\overline{x}) . \label{G-Bc-expr}
\een
For $G_i$ given by the first of these expressions, we can
see that function $\Gamma _i$ defined by (\ref{Gamma})
and the unique solution to equation (\ref{gG-G-eqn})
are given by
\begin{align}
\Gamma _i (x, \overline{x}) & = \biggl(
\frac{\alpha_i (n_i - 1)}{(r_i - \mu) n_i k_i} \biggr)
^{1/(1-\alpha_i)} \bigl( x \hsym (\overline{x})
\bigr) ^{1/(1-\alpha_i)}
\label{Gamma-expr} \\
\text{and} \quad
\gbG_i (\overline{x}) & = \biggl(
\frac{\alpha_i (n_i - 1)}{(r_i - \mu) n_i k_i} \biggr)
^{1/(1-\alpha_i)} \bigl( \overline{x} \hsym (\overline{x})
\bigr) ^{1/(1-\alpha_i)} . \label{gG-expr}
\end{align}
In particular, we note that $\gbG_i$ is as in
(\ref{eq-barX2}), while
\ben
0 < \Gamma_i (x, \overline{x}) \leq \Gamma_i
(\overline{x}, \overline{x}) = \gbG_i (\overline{x})
\quad \text{for all } 0 < x < \overline{x} .
\label{Gamma-Phi}
\een
Furthermore, Assumptions~\ref{A2} and~\ref{A4}
imply that $\{ \gbG_i , \ i \in I \} \in \bcal$ and,
given any $\varepsilon \in \mbox{} ]0, \varepsilon_i[$,
\be
0 < \gbG_i (\overline{x}) \leq \biggl( \widetilde{K}
(\varepsilon) \frac{\alpha_i (n_i - 1)}{(r_i - \mu) n_i k_i}
\biggr) ^{1/(1-\alpha_i)} \overline{x}
^{\frac{1 - \overline{\alpha}}{(1 - \overline{\alpha}
+ \overline{\alpha} \gamma) (1 - \alpha_i)}
+ \frac{\varepsilon}{1 - \alpha_i}}
\quad \text{for all } \overline{x} > \widetilde{x}
(\varepsilon) > 1 .
\ee
Combining this estimate with the identity
\be
\frac{(1 - \overline{\alpha}) \alpha_i}
{(1 - \overline{\alpha} + \overline{\alpha}
\gamma) (1 - \alpha_i)}
+ \frac{\alpha_i \varepsilon_i}{1 - \alpha_i}
= n_i - 1
\ee
and the inequality
\be
\frac{1 - \overline{\alpha}}
{(1 - \overline{\alpha} + \overline{\alpha}
\gamma) (1 - \alpha_i)}
+ \frac{\alpha_i \varepsilon_i}{1 - \alpha_i}
< n_i ,
\ee
which follow from Assumption~\ref{A2}, we can
see that there exists $\vartheta_i > 0$
such that
\begin{align}
x \gbG_i^{\alpha_i} (\overline{x}) & \leq
\biggl( \widetilde{K} (\vartheta_i)
\frac{\alpha_i (n_i - 1)}{(r_i - \mu) n_i k_i}
\biggr) ^{\alpha_i / (1-\alpha_i)} \overline{x}
^{n_i - \vartheta_i} \quad \text{for all }
\overline{x} > \widetilde{x} (\vartheta_i)
\label{Phi^a-est} \\
\text{and} \quad
\gbG_i (\overline{x}) & \leq \biggl( \widetilde{K}
(\vartheta_i) \frac{\alpha_i (n_i - 1)}
{(r_i - \mu) n_i k_i} \biggr) ^{1/(1-\alpha_i)}
\overline{x} ^{n_i - \vartheta_i}
\quad \text{for all } \overline{x} > \widetilde{x}
(\vartheta_i) . \label{Phi-est}
\end{align}
On the other hand, the solution to the ODE in
(\ref{G-Bc-expr}) is given by
\ben
B_i (c,\overline{x}) = \widetilde{f}_i (\overline{x}) +
\frac{k_i}{(n_i - 1) \bigl( n_i (1 - \alpha_i) - 1 \bigr)}
c^{-n_i (1 - \alpha_i) + 1} \overline{x} ^{\, -n_i}
\gbG _i^{n_i (1 - \alpha_i)} (\overline{x})
, \label{B-expr}
\een
where $\widetilde{f}_i$ is a function that we will
determine
later.

The requirement that $w$ should be continuous
is reflected by the identity
\be
\lim _{c \uparrow \gbG (\overline{x})} w(c, x, \overline{x})
= \lim _{c \downarrow \gbG (\overline{x})}
w(c, x, \overline{x}) ,
\ee
which gives rise to the expressions
\ben
\gbG _i^{\alpha_i} (\overline{x}) A_i (\overline{x})
= B_i \bigl( \gbG_i (\overline{x}), \overline{x} \bigr)
= \widetilde{f}_i  (\overline{x}) +
\frac{k_i}{(n_i - 1) \bigl( (1 - \alpha_i) n_i - 1 \bigr)}
\overline{x} ^{\, -n_i} \gbG _i (\overline{x})
. \label{A-expr}
\een

To determine the functions $A_i$ and $f_i $,
we first recall that, once the state process takes values
on the surface $\scal_i^1$, it should be optimal to be
kept there at all times by means of minimal control
effort exerted when the state process is on the boundary
line defined by $x = \overline{x}$ and  $\overline{x}
= \gbG _i (\overline{x})$.
This observation and the ``smooth-pasting condition''
of singular stochastic control suggest that $w$
should satisfy
\be
\lim _{c \downarrow \gbG_i (\overline{x}) , \,
x \uparrow \overline{x}}
w_c (c,x,\overline{x}) = k_i
\quad \text{and} \quad
\lim _{c \downarrow \gbG_i (\overline{x}) , \,
x \uparrow \overline{x}}
w_{cx} (c,x,\overline{x}) = 0 ,
\ee
which give rise to the system of equations
\begin{align}
\gbG _i^{\alpha_i -1} (\overline{x}) \Bigl( \alpha_i
A_i (\overline{x}) + \gbG_i (\overline{x}) f_i  \bigl(
\gbG_i (\overline{x}) \bigr) \Bigr) \overline{x}^{n_i}
+ \frac{\alpha_i}{r_i - \mu} \gbG _i^{\alpha_i-1}
(\overline{x}) \overline{x} \hsym (\overline{x})
& = k_i \label{FB1} \\
\text{and} \quad
n_i \gbG _i^{\alpha_i -1} (\overline{x}) \Bigl( \alpha_i
A_i (\overline{x}) + \gbG_i (\overline{x}) f_i  \bigl(
\gbG_i (\overline{x}) \bigr) \Bigr) \overline{x}^{n_i}
+ \frac{\alpha_i}{r_i - \mu} \gbG _i^{\alpha_i - 1}
(\overline{x}) \overline{x} \hsym (\overline{x})
& = 0 . \label{FB2}
\end{align}
The solution to this system is given by $\gbG_i
(\overline{x})$ as in (\ref{gG-expr}) and by
\ben
f_i  \bigl( \gbG_i (\overline{x}) \bigr) = - \frac{k_i}{n_i - 1}
\overline{x} ^{\, -n_i} \gbG _i^{-\alpha_i} (\overline{x})
- \alpha_i \gbG _i^{-1} (\overline{x}) A_i (\overline{x})
, \label{f-expr1}
\een
which is equivalent to
\ben
f_i  (y) = - \frac{k_i}{n_i - 1} y^{-\alpha_i} \bigl( \gbG
_i^\inv \bigr)^{-n_i} (y) - \alpha_i y^{-1}
A_i \bigl( \gbG _i^\inv (y) \bigr) , \label{f-expr2}
\een
where $\gbG _i^\inv$ is the inverse function of
$\gbG_i$.

To complete this part of the analysis, we note that
the requirement that $w$ should satisfy (\ref{vz-eq})
gives rise to the expression
\ben
\lim _{x \uparrow \overline{x}} w_{\overline{x}}
(c, x, \overline{x}) = c^{\alpha_i} \biggl( \Bigl(
(A_i) _-' (\overline{x}) - f_i  \bigl( \gbG_i
(\overline{x}) \bigr) (\gbG_i) _-' (\overline{x})
\Bigr) \overline{x} ^{n_i} + \frac{1}{r_i - \mu}
\overline{x} \hsym _-' (\overline{x}) \biggr)
= 0 . \label{FB3}
\een
Substituting the expression given by (\ref{f-expr1})
for $f_i $ in this equation and using the identities
\begin{align}
\frac{1}{r_i - \mu} \overline{x} \hsym (\overline{x})
\gbG _i^{\alpha_i} (\overline{x}) & =
\frac{n_i k_i}{\alpha_i (n_i - 1)} \gbG_i (\overline{x})
\label{gG-expr-impl} \\
\text{and} \quad
\frac{1}{r_i - \mu} \bigl( \overline{x} \hsym (\overline{x})
\bigr) _-' \gbG _i^{\alpha_i} (\overline{x}) & =
\frac{(1-\alpha_i) n_i k_i}{\alpha_i (n_i - 1)}
(\gbG_i) _-' (\overline{x}) , \nonumber
\end{align}
which follow from (\ref{gG-expr}), we obtain
\be
\bigl( \gbG _i^{\alpha_i} A_i \bigr) _-' (\overline{x})
= - \frac{(1-\alpha_i) k_i}{\alpha_i}
\overline{x} ^{\, -n_i} (\gbG _i) _-' (\overline{x})
- \frac{k_i}{\alpha_i (n_i - 1)} \bigl( \overline{x}
^{\, -n_i} \gbG _i (\overline{x}) \bigr) _-' .
\ee
In view of (\ref{Phi-est}) and the integration
by parts formula, we can see that both of
the integrals in the identity
\be
\overline{x}^{\, -n_i} \gbG _i (\overline{x})
= n \int _{\overline{x}}^\infty y^{\, -n_i - 1}
\gbG _i (y) \, \di y - \int _{\overline{x}}^\infty
y^{\, -n_i} (\gbG _i)_-' (y) \, \di y
\ee
converge.
It follows that
\begin{align}
A_i (\overline{x}) & = -
\frac{\bigl( (1-\alpha_i) (n_i-1) + 1 \bigr) k_i}
{\alpha_i (n_i - 1)} \overline{x}
^{\, -n_i} \gbG _i ^{1-\alpha_i} (\overline{x})
+ \frac{(1 - \alpha_i) n_i k_i}{\alpha_i}
\gbG _i ^{-\alpha_i} (\overline{x}) \int
_{\overline{x}}^\infty y^{-n_i-1} \gbG _i (y)
\, \di y \nonumber \\
& = - \frac{k_i}{\alpha_i (n_i - 1)} \overline{x}
^{\, -n_i} \gbG _i ^{1-\alpha_i} (\overline{x})
+ \frac{(1 - \alpha_i) k_i}{\alpha_i}
\gbG _i ^{-\alpha_i} (\overline{x}) \int
_{\overline{x}}^\infty y^{-n_i} (\gbG _i) _-'
(y) \, \di y . \label{A-exprF}
\end{align}

\begin{lem} \label{lem:v-sol}
Suppose that Assumptions~\ref{A1} and~\ref{A2}
hold true.
Also, let $\hsym$ be any function satisfying
Assumption~\ref{A4} and consider the free-boundary
functions  $G_i$ and $\gbG_i$ given by
(\ref{G-Bc-expr}) and (\ref{gG-expr}), as well as
the functions $A_i$, $f_i$ and $B_i$ given by
(\ref{B-expr}), (\ref{A-expr}), (\ref{f-expr2}) and
(\ref{A-exprF}).
The function $w$ given by (\ref{v-expr})
is a solution to the HJB equation (\ref{HJB}) with
boundary condition (\ref{HJB-BC}).
Furthermore, there exist constants $K > 0$
and $\vartheta_i > 0$ such that
\ben
0 < x \Gamma _i^{\alpha_i} (x, \overline{x}) \leq K
(1 + \overline{x} ^{n_i - \vartheta_i}) , \quad
0 < x \gbG _i^{\alpha_i} (\overline{x}) \leq K
(1 + \overline{x} ^{n_i - \vartheta_i}) ,
\label{G-est}
\een
and
\ben
0 < w (c, x, \overline{x}) \leq K \Bigl( 1 +
\overline{x} ^{n_i - \vartheta_i} +
(1 + \overline{x} ^{n_i - \vartheta_i}) \bigl(
c \gbG _i^{-1} (\overline{x}) \bigr) ^{\alpha_i}
+ \overline{x} ^{n_i - \vartheta_i} \bigl(
c^{-1} \gbG _i (\overline{x}) \bigr)
^{n_i (1-\alpha_i) - 1} \Bigr) \label{w-est}
\een
for all $(c, x, \overline{x}) \in \ccal$.
\end{lem}
\noindent {\bf Proof.}
We first note that the functions $\gbG$ and $G$
satisfy the requirements of Problem~\ref{FBP}
because $1-\alpha_i > 0$, $\hsym \in \vcal$
and $\{ \gbG_i , \ i \in I \} \in \bcal$ (see also
the discussion after (\ref{Phi-est})).
In what follows, we use $K > 0$ as a generic
constant.
The estimates in (\ref{G-est}) follow immediately
from the expressions of $\Gamma_i$ and $\gbG _i$
in (\ref{Gamma-expr}) and (\ref{gG-expr}), the
estimate (\ref{Phi^a-est}) and the facts that
$\gbG _i$ is continuous and $\lim
_{\overline{x} \downarrow 0} \gbG_i (\overline{x})
= 0$.
The estimate (\ref{Phi-est}) and the first
expression for $A_i$ in (\ref{A-exprF}) imply that
\be
\gbG_i^{\alpha_i} (\overline{x}) \bigl| A_i (\overline{x}) \bigr|
\leq K (\overline{x} ^{\, -n_i} + \overline{x}
^{\, -\vartheta_i})
\quad \text{and} \quad
\overline{x} \hsym (\overline{x}) \leq K \gbG
_i^{1 - \alpha_i} (\overline{x}) \leq \gbG_i^{-\alpha_i} (\overline{x}) \overline{x}
^{n_i - \vartheta_i} .
\ee
On the other hand, (\ref{A-expr}) implies that
\be
\bigl| \widetilde{f}_i  (\overline{x}) \bigr| \leq K
\Bigl( \gbG _i^{\alpha_i} (\overline{x}) \bigl|
A_i (\overline{x}) \bigr| + \overline{x} ^{\, -n_i}
\gbG _i (\overline{x}) \Bigr) \leq K (\overline{x}
^{\, -n_i} + \overline{x} ^{\, -\vartheta_i}) ,
\ee
while (\ref{B-expr}) yields
\be
\bigl| B_i (c, \overline{x}) \bigr| \leq K \Bigl(
\overline{x} ^{\, -n_i} + \overline{x}
^{\, -\vartheta_i} + c \overline{x} ^{\, - n_i}
\bigl( c^{-1} \gbG _i (\overline{x}) \bigr)
^{n_i (1-\alpha_i)} \Bigr) .
\ee
Using (\ref{f-expr1}) and the second expression
for $A_i$ in (\ref{A-exprF}), we obtain
\be
f_i  \bigl( \gbG_i (\overline{x}) \bigr) = -
(1 - \alpha_i) k_i \gbG _i ^{-1-\alpha_i}
(\overline{x}) \int _{\overline{x}}^\infty y^{-n_i}
(\gbG _i) _-' (y) \, \di y < 0 .
\ee
The upper bound in (\ref{w-est}) follows from
these estimates and the definition (\ref{v-expr}) of
$w$.

In view of the expression (\ref{v-expr}) of $w$
the function $w_{\overline{x}} \bigl( c, x, \cdot
\bigr) : \mbox{} ]x, \infty[ \mbox{} \rightarrow
\bbr$ is the difference of two convex functions
for all $c \geq 0$ and $x>0$.
By construction, the function $w _c$ is continuous
along the surface $\scal _i^1$.
The continuity of $w_x$ and $w_{xx}$
across the surface $\scal _i^1$ follows from the
observations that, given any $(c, x, \overline{x}) \in
\ccal_i \setminus \scal _i^1$,
\begin{align}
w_x (c, x, \overline{x}) & = \lim
_{c \downarrow \Gamma_i (x, \overline{x})}
w_x (c , x, \overline{x})
+ \Bigl( w_c \bigl( \Gamma_i (x, \overline{x}) ,
x, \overline{x} \bigr) - k_i \Bigr) (\Gamma_i) _x
(x, \overline{x}) \nonumber \\
& = \lim _{c \downarrow \Gamma_i (x, \overline{x})}
w_x (c , x, \overline{x})
= w_x \bigl( \Gamma_i (x, \overline{x}) , x,
\overline{x} \bigr) \label{v_x-inc} \\
\text{and} \quad
w_{xx} (c, x, \overline{x}) & = \lim
_{c \downarrow \Gamma_i (x, \overline{x})}
w_{xx} (c , x, \overline{x})
+ w_{cx} \bigl( \Gamma_i (x, \overline{x}) ,
x, \overline{x} \bigr) (\Gamma_i) _x (x, \overline{x})
\nonumber \\
& = \lim _{c \downarrow \Gamma_i (x, \overline{x})}
w_{xx} (c , x, \overline{x})
= w_{xx} \bigl( \Gamma_i (x, \overline{x}) , x,
\overline{x} \bigr) , \label{v_xx-inc}
\end{align}
where we have used the definition (\ref{v-expr}) of
$w$ as well as (\ref{v_c-G}) and (\ref{v_cx-G}).
The definition (\ref{v-expr}) of $w$ and the expression
(\ref{A-expr}) imply that $w_x$ and $w_{xx}$
are continuous across the surface $\scal _i^2$.
On the other hand, we use equation (\ref{gG-G-eqn})
as well as the expressions of $(B_i)_c$ and $f$ in
(\ref{G-Bc-expr}) and (\ref{f-expr1}) to obtain
\begin{align}
\lim _{c \downarrow \gbG_i (\overline{x})}
w_c (c, x, \overline{x}) & =
- \frac{k_i}{n_i - 1} \overline{x} ^{\, -n_i} x^{n_i}
+ \frac{\alpha_i}{r_i - \mu} \gbG _i^{\alpha_i-1}
(\overline{x}) x \hsym (\overline{x})
\nonumber \\
& = - \frac{k_i}{n_i - 1} G_i^{-n_i} \bigl( \gbG_i
(\overline{x}) , \overline{x} \bigr) x^{n_i}
+ \frac{\alpha_i}{r_i - \mu} \gbG _i^{\alpha_i-1}
(\overline{x}) x \hsym (\overline{x})
\nonumber \\
& = \lim _{c \uparrow \gbG_i (\overline{x})}
w_c (c, x, \overline{x}) . \nonumber
\end{align}
It follows that the function $w (\cdot, \cdot,
\overline{x})$ is $C^{1,2}$ in the interior of
$\ccal$ for all $\overline{x} > 0$.

By construction, we will prove that $w$ provides
a solution to Problem~\ref{FBP} if we show that
\begin{align}
w_c (c, x, \overline{x}) & < k_i \quad \text{for all }
0 < x < \overline{x} \text{ and } \gbG_i (\overline{x})
< c , \label{vc-ineq-u-p} \\
w_c (c, x, \overline{x}) & < k_i \quad \text{for all }
0 < x \leq \overline{x} \text{ and } \Gamma _i
(x, \overline{x}) < c \leq \gbG _i (\overline{x}) ,
\label{vc-ineq-d-p} \\
\text{and} \quad
\lscr_i w (c, x, \overline{x}) + c^{\alpha_i}
x \hsym (\overline{x}) & \leq 0 \quad \text{for all }
0 < x < \overline{x} \text{ and } c \leq \Gamma_i
(x, \overline{x}) . \label{v-ODE-ineq}
\end{align}
To this end, consider first any $(c, x, \overline{x})$
such that $0 < x \leq \overline{x} \text{ and }
\Gamma _i (x, \overline{x}) < c \leq \gbG _i
(\overline{x})$.
Using (\ref{v_cx-G}), we calculate
\begin{align}
w_{cx} (c, x, \overline{x}) & = n_i (B_i)_c
(c, \overline{x}) x^{n_i - 1} + \frac{\alpha_i}{r_i - \mu}
c^{\alpha_i - 1} \hsym (\overline{x})
\nonumber \\
& = \frac{\alpha_i}{r_i - \mu} c^{\alpha_i - 1}
\Biggl( 1 - \biggl( \frac{x}{G_i (c, \overline{x})}
\biggr)^{n_i - 1} \Biggr) \hsym (\overline{x}) > 0
\quad \text{for all } x \in \mbox{} \bigl] 0, G_i
(c,\overline{x}) \bigr[ . \nonumber
\end{align}
Combining this inequality with (\ref{v_c-G}),
we obtain (\ref{vc-ineq-d-p}).

Next, we consider any $(c, x, \overline{x})$
such that $0 < x \leq \overline{x} \text{ and }
\gbG _i (\overline{x}) < c$ and we use (\ref{FB3})
to derive
\be
w_{c \overline{x}} (c, x, \overline{x}) =
\frac{\alpha_i}{r_i - \mu} c^{\alpha_i - 1} x^{n_i}
\hsym _-' (\overline{x}) \bigl( x^{-n_i+1}
- \overline{x} ^{\, -n_i+1} \bigr) < 0 .
\ee
In view of this inequality and (\ref{vc-ineq-d-p}),
we can see that
\be
w_c (c, x, \overline{x}) = w_c \bigl( c, x, \gbG
_i^\inv (c) \bigr) - \int _{\overline{x}}^{\gbG _i^\inv (c)}
w_{c \overline{x}} (c,x,y) \, \di y < k_i .
\ee

To establish (\ref{v-ODE-ineq}), we use the
definition (\ref {v-expr}) of $w$, as well as
(\ref{v_x-inc}) and (\ref{v_xx-inc}), to obtain
\begin{align}
\lscr_i w & (c, x, \overline{x}) + c^{\alpha_i}
x \hsym (\overline{x}) \nonumber \\
& = \lscr_i w \bigl( \Gamma_i (x,\overline{x}) ,
z , x \bigr) + r_i k_i \bigl( \Gamma_i (x,\overline{x})
- c \bigr) + c^{\alpha_i} x \hsym (\overline{x})
\nonumber \\
& = - \Gamma _i^{\alpha_i} (x,\overline{x}) x
\hsym (\overline{x}) + r_i k_i \bigl( \Gamma_i
(x,\overline{x}) - c \bigr) + c^{\alpha_i} x
\hsym (\overline{x}) \nonumber \\
& = - \int _c^{\Gamma _i (x,\overline{x})} \bigl(
\alpha_i u^{\alpha_i - 1} x \hsym (\overline{x}) -r_i k_i
\bigr) \, \di u \quad \text{for all } c < \Gamma_i
(x, \overline{x}) . \label{L-ineq}
\end{align}
In view of (\ref{rbnm1}) and (\ref{rbnm2}), we can
see that the free-boundary point $G_i
(c,\overline{x})$ given by (\ref{G-Bc-expr}) is the
unique solution to the equation
\be
\int _0^{G_i (c,\overline{x})} y^{-m_i - 1} \bigl(
\alpha_i c^{\alpha_i - 1} y \hsym (\overline{x})
- r_i k_i \bigr) \, \di y = 0 .
\ee
Therefore, $\alpha_i c^{\alpha_i - 1} x
\hsym (\overline{x}) - r_i k_i > 0$ for all $x \geq
G_i (c,\overline{x})$, which implies that
\be
\alpha_i c^{\alpha_i - 1} x \hsym (\overline{x})
- r_i k_i > 0 \quad \text{for all } c \leq \Gamma_i
(x,\overline{x})
\ee
because $\Gamma_i (\cdot,\overline{x})$ is the
inverse of the strictly increasing function
$G_i (\cdot, \overline{x})$.
However, this conclusion and (\ref{L-ineq})
imply that $w$ satisfies (\ref{v-ODE-ineq}).

To prove the positivity of $w$ and complete the
proof, we first use (\ref{FB3}) to calculate
\begin{align}
\int _{\gbG_i (\overline{x})}^c f(y) \, \di y & =
\int _{\overline{x}}^{\gbG _i^\inv (c)} f \bigl(
\gbG_i (y) \bigr) (\gbG _i)_-' (y) \, \di y \nonumber \\
& = A_i \bigl( \gbG _i^\inv (c) \bigr) - A_i (\overline{x})
+ \frac{1}{r_i - \mu} \int _{\overline{x}}^{\gbG _i^\inv (c)}
y^{-n_i + 1} \hsym _-' (y) \, \di y . \nonumber
\end{align}
We next use integration by parts to obtain
\begin{align}
\int _{\overline{x}}^{\gbG _i^\inv (c)} y^{-n_i + 1}
\hsym _-' (y) \, \di y = \mbox{} & \bigl( \gbG _i^\inv
(c) \bigr) ^{-n_i + 1} \hsym \bigl( \gbG _i^\inv (c)
\bigr) - \overline{x} ^{\, -n_i + 1} \hsym (\overline{x})
\nonumber \\
& + (n_i - 1) \int _{\overline{x}}^{\gbG _i ^\inv (c)}
y^{-n_i} \hsym (y) \, \di y . \nonumber
\end{align}
On the other hand, we use (\ref{gG-expr-impl})
and the second expression for $A_i$ in
(\ref{A-exprF}) to obtain
\begin{align}
A_i (\overline{x}) + \frac{1}{r_i - \mu} \overline{x}
^{\, -n_i + 1} \hsym (\overline{x}) =
\mbox{} & \frac{k_i}{\alpha_i} \overline{x}
^{\, -n_i} \gbG _i ^{1-\alpha_i} (\overline{x})
\nonumber \\
& + \frac{(1 - \alpha_i) k_i}{\alpha_i}
\gbG _i ^{-\alpha_i} (\overline{x}) \int
_{\overline{x}}^\infty y^{-n_i} (\gbG _i) _-'
(y) \, \di y > 0 , \nonumber
\end{align}
which implies that
\be
A_i \bigl( \gbG _i^\inv (c) \bigr) + \frac{1}{r_i - \mu}
\bigl( \gbG _i^\inv (c) \bigr) ^{-n_i + 1} \hsym
\bigl( \gbG _i^\inv (c) \bigr) > 0 .
\ee
In view of these observations, we can see that,
given any $(c, x, \overline{x}) \in \ccal _i^3$,
\begin{align}
w(c, x & , \overline{x}) \nonumber \\
= \mbox{} & c^{\alpha_i} \biggl(
A_i \bigl( \gbG _i^\inv (c) \bigr) + \frac{1}{r_i - \mu}
\bigl( \gbG _i^\inv (c) \bigr) ^{-n_i + 1} \hsym
\bigl( \gbG _i^\inv (c) \bigr) + \frac{n_i - 1}{r_i - \mu}
\int _{\overline{x}}^{\gbG _i^\inv (c)} y^{-n}
\hsym (y) \, \di y \biggr) x^{n_i} \nonumber \\
& + \frac{1}{r_i - \mu} c^{\alpha_i} x^{n_i}
\bigl( x^{-n_i + 1} - \overline{x} ^{\, -n_i + 1}
\bigr) \hsym (\overline{x}) > 0 . \nonumber
\end{align}

Given any $(c, x, \overline{x}) \in \ccal _i^2$,
the fact that $(B_i)_c (c, \overline{x}) < 0$
(see (\ref{G-Bc-expr})) implies that
\begin{align}
w(c, x, \overline{x}) & \geq B_i \bigl( \gbG _i
(\overline{x}) , \overline{x} \bigr) x^{n_i} +
\frac{1}{r_i - \mu} c^{\alpha_i} x \hsym (\overline{x})
\nonumber \\
& = \gbG _i^{\alpha_i} (\overline{x}) A_i (\overline{x})
+ \frac{1}{r_i - \mu} c^{\alpha_i} x \hsym (\overline{x})
= w \bigl( \gbG _i (\overline{x}) , x, \overline{x}
\bigr) > 0 .
\end{align}
Finally, combining the observations that
\be
0 < w(c, x, \overline{x}) , \quad
0 < w_x (c, x, \overline{x})
\quad \text{and} \quad \lim
_{\substack{(c, x, \overline{x}) \in \ccal _i^2 \\\
(c, x, \overline{x}) \rightarrow (0,0,0)}}
w(c, x, \overline{x}) = 0
\ee
with (\ref{v_x-inc}), we can see that $w$ is
strictly positive in $\ccal_i^1$ as well.
\mbox{}\hfill$\Box$

\section{The solution to the control problem studied
in Sections~\ref{sec:singular-HJB} and~\ref{sec:HJB-sol}}
\label{sec:sol}

\begin{thm} \label{thm:IP-sol}
Consider the singular stochastic control problem
formulated in Section~\ref{sec:singular-HJB} and
suppose that Assumptions~\ref{A1} and~\ref{A2}
hold true.
Also, let $\hsym$ be any function satisfying
Assumption~\ref{A4}.
The value function $v^{(i)}$ of the optimisation problem
that the individual producer $i \in I$ faces, which is
defined by (\ref{IP-v}), identifies with the function
$w = w^{(i)}$ in Lemma~\ref{lem:v-sol}, namely,
\be
w (c_i, x, \overline{x}) = \sup _{C \in \acal_i}
J_{c_i,x,\overline{x}}^{(i)} (C) \quad \text{for all }
(c_i, x, \overline{x}) \in \ccal ,
\ee
where $J_{c_i,x,\overline{x}}^{(i)}$ is defined by
(\ref{Ji}).
Furthermore, the expansion process
\ben
C_i^\dagger (t) = c_i \vee \Bigl( \Gamma_i (\overline{X}_t,
\overline{x}) {\bf 1} _{\{ \overline{X}_t < \overline{x} \}}
+ \gbG_i \bigl( \overline{X}_t \bigr) {\bf 1}
_{\{ \overline{X}_t \geq \overline{x} \}} \Bigr) ,
\label{Cstar-maxX}
\een
where $\overline{X}$ is defined by (\ref{PXbar}), while
$\Gamma_i$ and $\gbG_i$ are given by
(\ref{Gamma-expr}) and (\ref{gG-expr}),
is admissible as well as optimal for the individual
producer $i \in I$.
\end{thm}
\noindent {\bf Proof.}
Fix any initial condition $(c_i, x, \overline{x})
\in \ccal$ and consider any expansion
process $C \in \acal_i$.
Using It\^{o}'s formula and the fact that the process
$\widetilde{X} = \overline{x} \vee \overline{X}$
is continuous and increases on the set
$\bigl\{ X_t = \widetilde{X}_t \bigr\}$, we obtain
\begin{align}
e^{-r_i T} w & \bigl( C(T-), X_T, \widetilde{X}_T
\bigr) \nonumber \\
= \mbox{} & w(c_i ,x,\overline{x}) + \int _0^T e^{-r_i t}
\lscr_i w \bigl( C(t), X_t, \widetilde{X}_t \bigr) \, \di t
+ \int _0^T e^{-r_i t} w_c \bigl( C(t), X_t,
\widetilde{X}_t \bigr) \, \di C_t^\cc \nonumber \\
& + \sum _{0 \leq t < T} e^{-r_i t} \Bigl( w \bigl(
C(t), X_t, \overline{X}_t \bigr) - w \bigl( C(t-),
X_t, \widetilde{X}_t \bigr) \Bigr) \nonumber \\
& + \int _0^T e^{-r_i t} w_{\overline{x}} \bigl(
C(t), \widetilde{X}_t, \widetilde{X}_t \bigr) \,
\di \widetilde{X}_t + M_T , \nonumber
\end{align}
where
\ben
M_T = \int _0^T e^{-r_i t} \sigma X_t w_x \bigl(
C(t), X_t, \widetilde{X}_t \bigr) \, \di W_t . \label{M}
\een
In view of this identity, the positivity of $w$ and
the fact that $w$ satisfies the HJB equation
(\ref{HJB}) with boundary condition (\ref{HJB-BC}),
we can see that
\begin{align}
\int _0^T e^{-r_i t} C^{\alpha_i} (t) X_t
\hsym (\widetilde{X}_t) & \, \di t - k_i \int _{[0,T[}
e^{-r_i t} \, \di C(t)
\nonumber \\
= \mbox{} & w(c_i, x,\overline{x}) - e^{-r_i T}
w \bigl( C(T-), X_T, \widetilde{X}_T \bigr)
\nonumber \\
& + \int _0^T e^{-r_i t} \Bigl( \lscr_i w \bigl(
C(t), X_t, \widetilde{X}_t \bigr) + C^{\alpha_i}
(t) X_t \hsym (\widetilde{X}_t) \Bigr) \, \di t
\nonumber \\
& + \int _0^T e^{-r_i t} \Bigl( w_c \bigl(
C(t), X_t, \widetilde{X}_t \bigr) - k_i \Bigr)
\, \di C^\cc (t) \nonumber \\
& + \sum _{0 \leq t < T} e^{-r_i t} \int
_0^{\Delta C(t)} \Bigl( w_c \bigl( C(t-) + y ,
X_t, \widetilde{X}_t \bigr) - k_i \Bigr) \, \di y
+ M_T \nonumber \\
\leq \mbox{} & w(c_i ,x, \overline{x}) - e^{-r_i T}
w\bigl( C(T-), X_T, \widetilde{X}_T \bigr)
+ M_T . \nonumber
\end{align}
Taking expectations, we are faced with the
inequality
\begin{align}
\EXP \biggl[ \int _0^{T \wedge \tau_n} e^{-r_i t}
C^{\alpha_i} (t) & X_t \hsym (\widetilde{X}_t)
\, \di t - k_i \int _{[0, T \wedge \tau_n[} e^{-r_i t}
\, \di C(t) \biggr] \nonumber \\
& \leq w(c_i, x, \overline{x}) - \EXP \Bigl[
e^{-r_i (T \wedge \tau_n)} w \bigl(
C((T \wedge \tau_n)-), X_{T \wedge \tau_n},
\widetilde{X}_{T \wedge \tau_n} \bigr) \Bigr]
, \label{VT-ineq}
\end{align}
where $(\tau_n)$ is a localising sequence for
the local martingale $M$.
Recalling the admissibility condition (\ref{Ac-IC}),
the fact that $\hsym$ is bounded and the positivity
of $w$, we can pass to the limit as
$n , T \rightarrow \infty$ using the monotone
convergence theorem to obtain $J_{c_i,x,
\overline{x}}^{(i)} (C) \leq w (c_i, x, \overline{x})$.
It follows that $v^{(i)} (c_i, x, \overline{x})
\leq w (c_i, x, \overline{x})$ because $C$
has been an arbitrary admissible expansion process.

To establish the reverse inequality, we first
note that the expansion process given
by (\ref{Cstar-maxX}) is such that
\be
{C_i^\dagger}^{\alpha_i} (t) \leq c_i^{\alpha_i}
+ \gbG _i^{\alpha_i} (\overline{x})
+ \gbG _i^{\alpha_i} \bigl( \overline{X}_t \bigr) .
\ee
Combining this inequality with the estimates
in (\ref{G-est}), the assumption that $n_i > 1$
and Lemma~1 in \cite{MZ07}, we can see that
\be
\EXP \biggl[ \int _0^\infty e^{-r_i t}
{C_i^\dagger}^{\alpha_i} (t) X_t \, \di t \biggr]
< \infty .
\ee
Therefore, $C^\dagger$ is admissible because
it satisfies the admissibility condition
(\ref{Ac-IC}).
This capacity schedule is such that the
(\ref{VT-ineq}) holds with equality, namely,
\begin{align}
\EXP \biggl[ \int _0^{T \wedge \tau_n} e^{-r_i t}
{C_i^\dagger}^{\alpha_i} (t) & X_t \hsym
(\widetilde{X}_t) \, \di t - k_i \int
_{[0, T \wedge \tau_n[} e^{-r_i t} \, \di
{C_i^\dagger} (t) \biggr] \nonumber \\
& = w(c_i ,x, \overline{x}) - \EXP \Bigl[
e^{-r_i (T \wedge \tau_n)} w \bigl( C_i^\dagger
(T \wedge \tau_n)-), X_{T \wedge \tau_n},
\widetilde{X} _{T \wedge \tau_n} \bigr) \Bigr]
. \label{VT-eq}
\end{align}

Using the expressions (\ref{Gamma-expr})
and (\ref{gG-expr}) of $\Gamma_i$ and
$\gbG _i$, as well as the estimates
in (\ref{w-est}), we can see that
\begin{align}
0 < e^{-r_i T} w \bigl( C_i^\dagger (T) , X_T,
\widetilde{X}_T \bigr) & \leq K e^{-r_i T} \Bigl(
1 + w \bigl( \gbG_i (\overline{X}_t) , X_T ,
\overline{X}_T \bigr) \Bigr) \nonumber \\
& \leq K e^{-r_i T} \bigl( 1 +
\overline{X} _T^{n_i - \vartheta_i} \bigr) ,
\label{VT-est}
\end{align}
where $K>0$ is a generic constant.
Therefore, we can pass to the limit as
$n \rightarrow \infty$ in (\ref{VT-eq}) using
the monotone and the dominated convergence
theorems to obtain
\begin{align}
\EXP \biggl[ \int _0^T e^{-r_i t}
{C_i^\dagger}^{\alpha_i} (t) X_t \hsym
(\widetilde{X}_t) \, \di t - k_i \int _{[0, T[} &
e^{-r_i t} \, \di C_i^\dagger (t) \biggr] \nonumber \\
& = w(c_i, x, \overline{x}) - \EXP \Bigl[ e^{-r_i T}
w \bigl( C_i^\dagger (T-), X_T, \widetilde{X}_T
\bigr) \Bigr] . \nonumber
\end{align}
Furthermore, we can use the admissibility of
the expansion process $C_i^\dagger$,
(\ref{VT-est}), Lemma~1 in~\cite{MZ07} and the
monotone and dominated convergence theorems
to pass to the limit as $T \uparrow \infty$
and obtain $J_{c_i,x, \overline{x}}^{(i)}
\bigl( C_i^\dagger \bigr) = w (c_i, x, \overline{x})$,
which implies that $v^{(i)} (c_i ,x, \overline{x})
\geq w (c_i, x, \overline{x})$ as well as the
optimality of $C_i^\dagger$.
\mbox{}\hfill$\Box$

\section{Proof of Theorem~\ref{thm:main}}
\label{sec:main-proof}

The following observations prove that the
expressions (\ref{eq.PbarP*}) and (\ref{eq.C*}),
or equivalently, (\ref{eq-barX1}) and (\ref{eq-barX2}),
provide a production equilibrium in the sense of
Definition~\ref{de:equil}.

\begin{description}
\item[1.]
Consider any family $\bgbG = \{ \gbG_i , \ i \in I \}
\in \bcal$, where $\bcal$ is introduced by
Definition~\ref{S-Bcal}, and suppose that
the function $\hsym = \iscr [\bgbG]$
satisfies the requirements of
Assumption~\ref{A4}.

\item[2.]
With the possible exception of the producer indexed
by $i \in I$, suppose that every other producer, say
$j \in I \setminus \{ i \}$, adopts the capacity
expansion strategy
$C_j (t) = c_j \vee \gbG_j (\overline{X}_t)$.
In the competitive setting that we have considered,
this assumption and the clearing condition
(\ref{Peq1}) give rise to the price
process 
\ben
P = \bigl( X \hsym (\overline{X}) \bigr) ^{1 / (1+\beta)}
. \label{P-finpr}
\een

\item[3.]
Producer~$i$, who has been singled
out, is faced with the task of maximising the
performance index
\be
\widetilde{J}_{c_i}^{(i)} (C_i \mid P) =
J_{c_i, \overline{x}, \overline{x}}^{(i)} (C_i)
\ee
over all $C_i \in \acal$,
where $\widetilde{J}_{c_i}^{(i)} (C_i \mid P)$
and $J_{c_i, x, \overline{x}}^{(i)} (C_i)$
are the performance criteria defined by
(\ref{Ji-tilde}) and (\ref{Ji}).

\item[4.]
Theorem~\ref{thm:IP-sol} implies that the expansion
strategy
\begin{align}
C_i^\dagger (t) & = \biggl(
\frac{\alpha_i (n_i - 1)}{(r_i - \mu) n_i k_i} \biggr)
^{1/(1-\alpha_i)} \bigl( \overline{x} \hsym (\overline{x})
\bigr) ^{1/(1-\alpha_i)} \nonumber \\
& = c_i \vee \Usym_i^\star \Bigl( \bigl( \overline{x}
\hsym (\overline{x}) \bigr) ^{1/(1+\beta)} \Bigr)
\stackrel{(\ref{P-finpr})}{=} c_i \vee \Usym_i^\star
(\overline{P}_t) , \nonumber
\end{align}
where
\ben
\Usym_i^\star \bigl( \overline{p} \bigr) = \biggl(
\frac{\alpha_i (n_i - 1)}{(r_i - \mu) n_i k_i} \biggr)
^{1/(1-\alpha_i)} \overline{p} ^{(1+\beta) / (1-\alpha_i)} ,
\label{Psii*-finpr}
\een
is such that (\ref{Cistar-eq-de}) holds true.
In particular, this expansion strategy is optimal for
producer~$i$.

\item[5.]
The family $\bUsym^\star = \{ \Usym_i^\star, \ i \in I \}$
defined by (\ref{Psii*-finpr}) belongs to $\bcal$.
By Theorem~\ref{thm:Phi-Psi-eq}, the solution
to the equation
$P^{1+\beta} = X \iscr [\bUsym^\star] (\overline{P})$
is given by
$P^{1+\beta} = X \iscr [\bgbG^\star] (\overline{X})$,
where
\be
\gbG _i^\star (\overline{x}) = \Usym _i^\star \Bigl(
\bigl( \overline{x} \hsym^\star (\overline{x}) \bigr)
^{1/(1+\beta)} \Bigr)
\ee
and the function $\hsym^\star$ satisfies
$\hsym^\star = \iscr [\bgbG^\star]$ and
possesses the properties listed
in (\ref{phi-props}).

\item[6.]
The function $\hsym^\star$ all of the
requirements of Assumption~\ref{A4}, by virtue of
Lemma~\ref{lem:est} and (\ref{Psii*-finpr}).

\item[7.]
In light of the observations above, the price
process
\ben
P^\star = \bigl( X \hsym^\star (\overline{X}) \bigr)
^{1 / (1+\beta)}  = \bigl( X \usym^\star
(\overline{P}) \bigr)^{1 / (1+\beta)} ,
\label{P*-finpr}
\een
where $\usym^\star = \iscr [\bUsym^\star]$, and
the expansion strategies
\be
C_i^\star (t) = c_i \vee \Usym_i^\star \bigl(
\overline{X}_t \, \hsym (\overline{X}_t) \bigr)
^{1 / (1-\alpha_i)}  = c_i \vee \Usym_i^\star
(\overline{P} _t^\star)
\ee
comprise a competitive production equilibrium.

\end{description}

\noindent
Finally, (\ref{dP-eq}) and (\ref{P0-eq}) follow
immediately from (\ref{dX}), (\ref{P*-finpr}) and
an application of It\^{o}'s formula.
\mbox{}\hfill$\Box$

\small

\end{document}